\documentclass{article}
\usepackage{graphicx} 
\usepackage{tikz}
\usepackage{amsmath}
\usepackage{amssymb}

\usepackage{float}
\usepackage{wrapfig}
\usepackage{amsthm,mathtools}
\usepackage{thmtools}
\usepackage{setspace}
\usepackage[numbers]{natbib}
\usepackage{mathrsfs}
\usepackage{multicol}

\usepackage[margin=1.03in]{geometry} 

\usetikzlibrary{intersections,through}

\newcommand{\E}{\mathbb{E}}
\newcommand{\W}{\mathbf{W}}

\newcommand{\sgn}{\operatorname{sgn}}

\newtheorem{definition}{\hspace{2em}Definition}
\newtheorem{theorem}{\hspace{2em}Theorem}

\newtheorem{conjecture}{\hspace{2em}Conjecture}
\newtheorem{observation}{\hspace{2em}Observation}

\newtheorem{proposition}{\hspace{2em}Proposition}
\newtheorem{lemma}{\hspace{2em}Lemma}
\newtheorem{corollary}{\hspace{2em}Corollary}

\newtheorem{geometryproposition}{\hspace{2em}Proposition}[section]

\begin{document} 
	\title{A Two-regime Khintchine Inequality and an Improved Bound on the Degree-1 Fourier Weight for Linear Threshold Functions} 
	\author{Xuan Fang\footnote{26110180010@m.fudan.edu.cn} \quad and \quad Tianyu Wang\footnote{wangtianyu@fudan.edu.cn}}
	\date{} 
	\maketitle 
	
	\begin{abstract}
The Khintchine inequality provides a lower bound on the expected absolute value of a weighted sum of independent Rademacher random variables. In the classical setting, when the coefficient vector has unit Euclidean norm, this lower bound is a constant, with equality attained only by a simple family of extremal configurations. A refinement due to De, Diakonikolas, and Servedio, referred to as the \emph{linear Khintchine inequality} in this work, strengthens this estimate by establishing a lower bound that depends linearly on the distance of the coefficient vector from the extremal set.

In this paper, we present a refined analysis of this dependence on the weight vector. Our results reveal a phase transition in the rate of improvement: when the dimension exceeds six, the lower bound undergoes an abrupt change as the weight vector deviates from the minimizer. Additionally, we improve the slope constant in linear Khintchine inequality. As a consequence, we establish an improved lower bound on the degree-1 Fourier weight for linear threshold functions $\mathbf W^{\leq1}[f]\geq0.55100929$, marking
progress towards a conjecture of Benjamini, Kalai, and Schramm.


\end{abstract}

\section{Introduction}
	
	Khintchine-type inequalities are fundamental tools for estimating random signed sums. Let $\epsilon_1,\ldots,\epsilon_n$ be independent Rademacher random variables and let $\omega=(\omega_1,\ldots,\omega_n)\in\mathbb{R}^n$.  The quantity
	\[
	K_n(\omega)=\E\left|\sum_{i=1}^n \omega_i\epsilon_i\right|
	\]
	measures the average size of the linear form $\omega\cdot x$ on the Hamming cube.  The classical Khintchine inequality, going back to Khintchine \cite{khintchine1923dyadische}, asserts the equivalence of the $L_p$ and $L_2$ norms of such sums; the sharp constants in this family were subsequently studied by Szarek, Haagerup, and others \cite{szarek1976best,haagerup1981best,latala1994best,nayar2012khintchine}. 
    
	For the endpoint of primary relevance to the present work, the sharp lower bound is 
    \[
	K_n(\omega)\geq \frac{1}{\sqrt2}\|\omega\|_2.
	\]
	After the normalization $\|\omega\|_2=1$, the best possible value is therefore $1/\sqrt2$.  The equality cases are highly rigid: equality is attained only when, up to signs and permutations, $\omega$ has two nonzero coordinates, both equal to $1/\sqrt2$.  This rigidity naturally leads to a stability question: if $\omega$ is separated from the extremal set, how would $K_n(\omega)$ quantitatively deviate from $1/\sqrt2$? Previously, this stability question was studied by De, Diakonikolas, and Servedio \cite{de2016robust} and by Eskenazis, Nayar, and Tkocz \cite{eskenazis2024resilience}.
	
	A motivation for this question comes from the Fourier analysis of Boolean functions.  For a Boolean function $f:\{-1,1\}^n\to\{-1,1\}$, write
	\[
	\mathbf{W}^{\leq 1}[f]=\widehat f(\emptyset)^2+
	\sum_{i=1}^n\widehat f(i)^2
	\]
	for its Fourier weight on levels $0$ and $1$.  If $f(x)=\sgn(\ell(x))$ is a linear threshold function (LTF) with linear function $\ell$, then the degree-$0$ and degree-$1$ Fourier coefficients, also called Chow parameters, encode strong structural information about $f$ \cite{odonnell2011chow,de2014nearly}.  Benjamini, Kalai, and Schramm first posed the conjecture \cite[Remark~3.7]{benjamini1999noise}, subsequently recorded in O'Donnell's book \cite[Conjecture~5.3]{odonnell2014analysis}, that every LTF $f$ satisfies 
	\begin{align} 
	    \mathbf{W}^{\leq1}[f]\geq \frac{2}{\pi}. \label{eq:bks-conj} 
	\end{align} 
	The constant in \eqref{eq:bks-conj} is asymptotically attained by majority functions; see also \cite{benjamini1999noise,gotsman1994spectral}. The standard proof of the weaker bound $\mathbf{W}^{\leq1}[f]\geq1/2$ uses Khintchine's inequality through the estimate 
	\[ 
	\sqrt{\mathbf{W}^{\leq1}[f]}
	\geq \frac{\E|\ell(x)|}{\|\ell\|_2}
	\] 
	(after the usual normalization, and with the constant term handled separately).
    
    In the works toward this conjecture on the degree-$1$ Fourier weight, the linear Khintchine inequalities emerged as a crucial ingredient: they separate the near-extremal two-coordinate configurations from genuinely higher-dimensional coefficient vectors and thus provide the additional quantitative information that the classical inequality alone does not capture. This is a principal reason for studying the sharp stability profile of $K_n$ in the present paper.

	The first systematic study on the linear Khintchine inequality was developed by De, Diakonikolas, and Servedio \cite{de2016robust}.  They proved that if $\tau := \tau(\omega)$ denotes the Euclidean distance from a normalized coefficient vector to the extremal set of the classical Khintchine inequality, then
	\begin{align}
	    K_n(\omega)\geq \frac1{\sqrt2}+c\tau (\omega) \label{eq:khintchine-dds}
	\end{align} 
	for a universal constant $c>0$. We refer to inequalities of this type as linear or robust Khintchine inequalities. This result yielded a strict improvement over the $1/2$ lower bound for the degree-$1$ Fourier weight of LTFs and also led to algorithms for computing several optimal constants in Fourier analysis and high-dimensional geometry.  Eskenazis, Nayar, and Tkocz subsequently extracted explicit constants from a quantitative stability proof \cite{eskenazis2024resilience}.  

    More recently, Rademacher sums and Khintchine-type inequalities have been studied from several complementary viewpoints.  Keller and Klein resolved Tomaszewski's anti-concentration conjecture using local concentration inequalities and refined Berry--Esseen estimates \cite{keller2022proof}; Kalarickal, Rotunno, Singh, and Tkocz established convexity properties for Rademacher sums on the ordered chamber \cite{kalarickal2025convexity}; and a growing body of work has developed stability or deficit versions of sharp Khintchine-type inequalities for Rademacher and related distributions \cite{eskenazis2024distributional,jakimiuk2026stability,jakimiuk2026spheres,chavez2026stability}.  These developments indicate that the geometry of near-extremizers is a useful and active perspective. 

\subsection{Our results}
As a standard convention, we work in the proper chamber
\[
\mathfrak C=\{\omega:\|\omega\|_2=1,\ \omega_1\geq\cdots\geq\omega_n\geq0\},
\qquad \omega^*=(1/\sqrt2,1/\sqrt2,0,\ldots,0),
\]
and put $\mathcal C_n(\tau)=\{\omega\in\mathfrak C:
\|\omega-\omega^*\|_2=\tau\}$ and
$M_n(\tau)=\min_{\omega\in\mathcal C_n(\tau)}K_n(\omega)$.
Signs and permutations reduce the general case to this chamber.

\textbf{A two-regime Khintchine inequality and sharp profiles.} Our main result is a two-regime Khintchine inequality. 
We show that for $n \ge 6$, the quantity $K_n (\omega)$ undergoes an abrupt phase transition at 
\begin{align*}
    \tau = \| \omega - \omega^* \|_2 = d_{\max} (4) := \sqrt{2 - \sqrt{2}}, 
\end{align*}
where $\omega^* = \left( \frac{1}{\sqrt{2}} , \frac{1}{\sqrt{2}} , 0, \cdots, 0 \right)$ is the unique minimizer in $ \mathfrak C$. Concretely, we establish the following inequality
\begin{align}
    K_n (\omega) \ge 
    \begin{cases} 
        \frac{1}{\sqrt{2}} + \left( \sqrt{1 - \frac{\tau^2}{4}} - \tau  \right) \frac{\tau}{2 \sqrt{2}}, &\quad \text{for $\omega \in \mathfrak C$ }, \\ 
        \frac{3 (2 - \tau^2) }{4\sqrt{2}} + \frac{3}{8\sqrt{2}} \sqrt{1 - \frac{(2 - \tau^2)^2}{2} } , &\quad \text{for $ \omega \in \mathfrak C $ and $\tau \in \left[ \sqrt{2 - \sqrt{2} }, \sqrt{ 2 - \frac{\sqrt{966}}{22}} \right] $} 
    \end{cases}, \label{eq:two-regime-intro}
\end{align}
and identify the sharp profiles at which the strict equality is obtained. Equality in the first inequality of \eqref{eq:two-regime-intro} is attained by a profile with four nonzero coordinates, which we call the four-coordinate profile.
Likewise, the profile at which the second inequality takes strict equality in \eqref{eq:two-regime-intro} is referred to as the six-coordinate profile.

\begin{figure}[t]
\centering
\includegraphics[width=.76\linewidth]{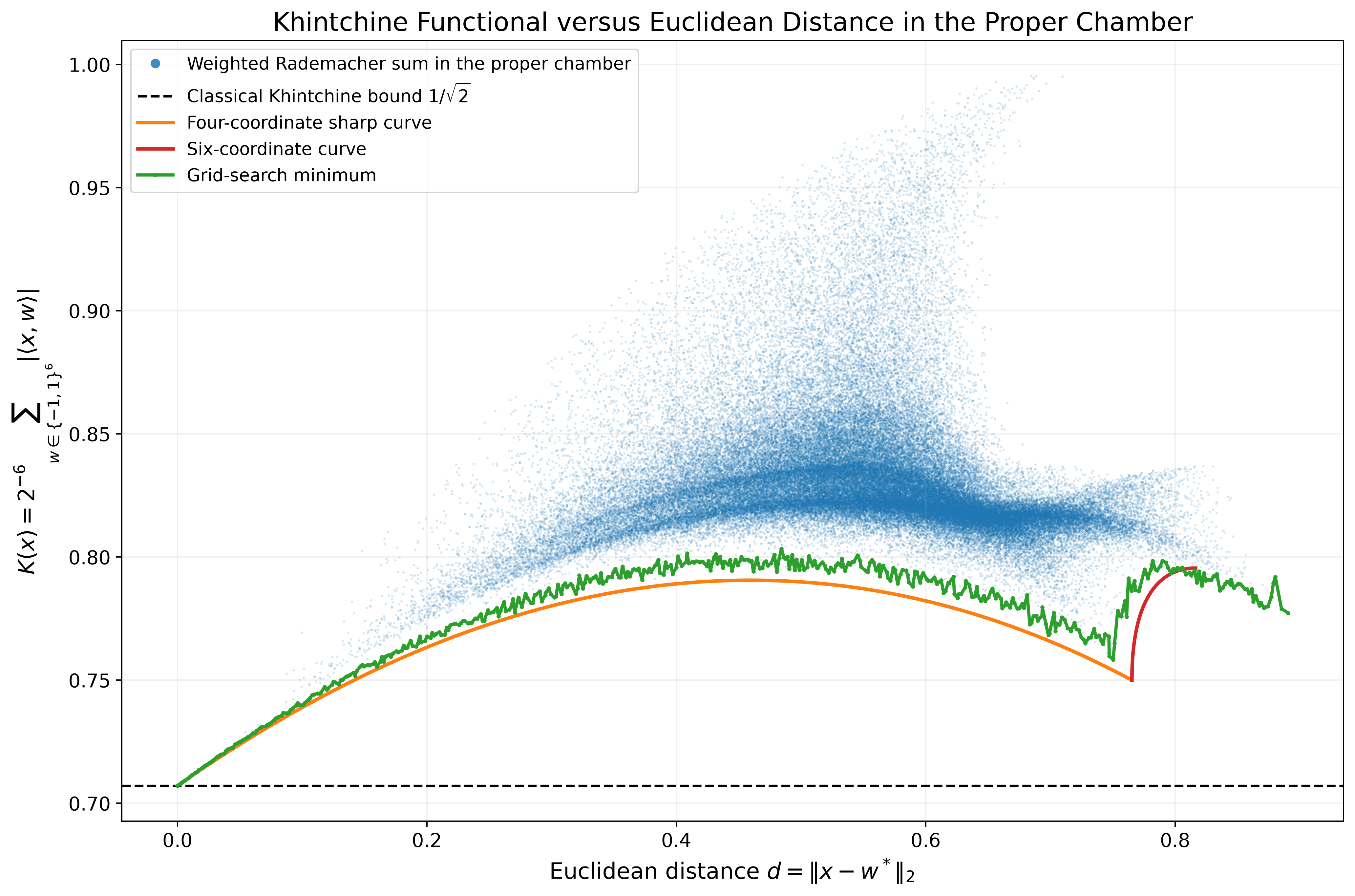}
\caption{Illustrations of the two-regime Khintchine inequality in dimension six.
The blue dots are $100{,}000$ points sampled uniformly from
$\mathbb S^5$ and folded into the proper chamber. The grid-search
curve is numerical; the four- and six-coordinate curves are the candidate
profiles. Theorems~\ref{thm:four} and \ref{th:six-local} certify the
four-coordinate interval and an explicit initial six-coordinate interval,
respectively.}
\label{fig:phase-transition}
\end{figure}

The inequality \eqref{eq:two-regime-intro} and the sharp profiles show an phase transition at $ \tau = \sqrt{2 - \sqrt{2}} $. This phase transition phenomenon is illustrated in Figure \ref{fig:phase-transition}.



\textbf{A linear Khintchine inequality with slope $1/20$.}
We prove a linear Khintchine inequality with the best known slope constant, to the best of our knowledge. This inequality states 
\begin{equation}\label{eq:improved-khintchine-linear}
 K_n(\omega)\geq\frac1{\sqrt2}+\frac1{20}\tau(\omega)
 \qquad(\|\omega\|_2=1).
\end{equation}
The proof combines Haagerup's refined inequality, a sharper Wallis-ratio estimate, and a polynomial minorant that uses the second, fourth, and sixth moments of a Rademacher tail. The four-coordinate extremizer gives the upper bound
\[
 c_{\mathrm{opt}}\leq
 \frac{3/4-1/\sqrt2}{\sqrt{2-\sqrt2}}=0.0560426911\ldots;
\]
our result gives $c_{\mathrm{opt}}\geq0.05$. This improves the explicit
coefficient $8\cdot10^{-5}$ in \cite{eskenazis2024resilience}.
The optimal coefficient remains conjectural.

\textbf{Low-degree Fourier weight of LTFs.}
Theorem~\ref{th:explicit-ltf} proves that every Boolean LTF $f$ satisfies
\[
 \mathbf W^{\leq1}[f]\geq\left(\frac{7423}{10000}\right)^2
 =0.55100929>0.55.
\]
This explicit improvement over $1/2$ follows from common coefficient
estimates, tie-free homogenization, and a conditioned Chow argument.
It does not require the robust theorem as an input. The conjectured
$2/\pi$ lower bound remains open.

\subsection{Organization}
Section~\ref{sec:profiles} proves the four- and six-coordinate profile
theorems. Section~\ref{sec:robust} proves the linear Khintchine inequality,
and Section~\ref{section:LTF} gives the Fourier-weight application.
Only the proofs of the four main theorems appear in the main text.
Supporting proofs are collected in Appendix~\ref{app:supporting}.
Appendix~\ref{app:phase-geometry} explains the geometry behind the proposed
higher phase transitions, while Appendix~\ref{app:additional} contains
additional local-minimizer and regularity properties.

\section{The two-regime Khintchine inequality and the sharp minimizing profiles}\label{sec:profiles} 

Let us start by reiterating notations and conventions that will be used throughout. For a unit vector $\omega\in\mathbb R^n$ with $n\geq2$, define
\begin{equation}\label{eq:proper}
	\tau=\tau(\omega):=\min_{1\leq i<j\leq n}\|\omega\pm\frac{1}{\sqrt{2}}e_i\pm\frac{1}{\sqrt{2}}e_j\|_2,
\end{equation}
where $e_i$ is the $i$-th coordinate vector and the minimum includes all choices of signs. Owing to the invariance of \(K_n(\omega)\) and \(\tau(\omega)\) under sign changes and coordinate permutations, there is no loss of generality in assuming that \(\omega\) is proper, i.e.,
\[
\omega \in \mathfrak{C} := \{ \omega : \omega_1 \ge \omega_2 \ge \cdots \ge \omega_n \ge 0, \quad \|\omega\|_2 = 1 \} .
\]
Within this chamber, the unique minimizer of the classical Khintchine inequality is
\[
\omega^* = \left(\frac{1}{\sqrt{2}}, \frac{1}{\sqrt{2}}, 0, \ldots, 0\right),
\]
and the distance to optimum is
\begin{align}
    \tau (\omega) = \|\omega - \omega^*\|_2. \label{eq:def-tau}
\end{align}
For each admissible \(\tau\), we denote by \(\mathcal{C}_n(\tau)\) the set of proper unit vectors at distance \(\tau\) from \(\omega^*\),
\begin{align*}
    \mathcal{C}_n(\tau) := \{ \omega : \| \omega - \omega^* \|_2 = \tau, \omega \in \mathfrak C \}. 
\end{align*}

\subsection{The four-coordinate profile}
Our first theorem establishes the first inequality in \eqref{eq:two-regime-intro}, and identifies the profile at which the minimum is obtained. 

\begin{theorem}\label{thm:four}
    For every $\omega\in\mathcal{C}_n(\tau)$,
    \[
    K_n(\omega)
    \geq \frac{1}{\sqrt2}
    +\left(\sqrt{1-\frac{\tau^2}{4}}-\tau\right)\frac{\tau}{2\sqrt2}.
    \]
    Moreover, when $n\geq4$, the estimate is optimal for $0\leq \tau\leq\sqrt{2-\sqrt2}$.  In that range equality holds if and only if
    \[
    \omega=\omega_\tau:=(r,r,s,s,0,\ldots,0),
    \]
    where
    \[
    r=\frac{1-\tau^2/2}{\sqrt2},
    \qquad
    s=\frac{\tau}{\sqrt2}\sqrt{1-\frac{\tau^2}{4}}.
    \]
\end{theorem}



	\begin{proof}[Proof of Theorem \ref{thm:four}]
Put $a=1/\sqrt2-\omega_1$, $b=1/\sqrt2-\omega_2$,
$c^2=\sum_{i\geq3}\omega_i^2$, and $S=\sum_{i\geq3}\omega_i\epsilon_i$.
Normalization and the distance constraint give
\begin{equation}\label{eq:tau-basic}
 \tau^2=a^2+b^2+c^2=\sqrt2(a+b).
\end{equation}
For $q=\omega_1+\omega_2$ and $p=\omega_1-\omega_2\geq0$,
averaging the first two signs gives
\[
 K_n(\omega)=\frac12\E\max\{q,|S|\}
             +\frac12\E\max\{p,|S|\}.
\]
The first expectation is at least $q$. By symmetry, the second is
$\E|p\epsilon_0+S|$, where $\epsilon_0$ is an independent Rademacher
variable. The sharp Khintchine inequality gives
\[
 \E\max\{p,|S|\}\geq\frac1{\sqrt2}
       \sqrt{p^2+c^2}.
\]
Using \eqref{eq:tau-basic}, we have
\[
 p^2+c^2=\tau^2-\frac{\tau^4}{2}+a^2+b^2
 \geq\tau^2-\frac{\tau^4}{2}+\frac{(a+b)^2}{2}
 =\tau^2-\frac{\tau^4}{4}.
\]
Since $q=\sqrt2-\tau^2/\sqrt2$, the claimed inequality follows.

For $\tau>0$, equality requires $a=b$ and equality in the tail
Khintchine inequality. Thus precisely two tail coefficients are nonzero
and equal, giving $(r,r,s,s,0,\ldots,0)$ with the stated parameters.
At $\tau=0$ the vector is $\omega^*$. Conversely, the candidate has
$\|\omega\|_2=1$ and the prescribed distance, and $r\geq s$ exactly
when $\tau\leq\sqrt{2-\sqrt2}$. In that range $|S|\leq2s\leq2r=q$,
so all inequalities above are equalities. This also proves uniqueness.
\end{proof}

\begin{proposition}\label{prop:dmax}
        Let
		\[
		d_{\max}(n)=\max_{\omega\in\mathbb{R}^n,\ \|\omega\|_2=1}\tau(\omega).
		\]
		Equivalently, after reduction to the proper chamber, $d_{\max}(n)$ is the maximum of $\|\omega-\omega^*\|_2$ over proper unit vectors in $\mathbb{R}^n$.
		For $n\geq2$,
		\[
		d_{\max}(n)=
		\begin{cases}
			\sqrt{2-\sqrt2}, & 2\leq n\leq4,\\[2mm]
			\sqrt{2-2\sqrt{2/n}}, & n\geq5.
		\end{cases}
		\]
		The two formulas agree at $n=4$.
	\end{proposition}
	
	The proof is given in Appendix~\ref{app:proof-prop-dmax}.

	Theorem \ref{thm:four} suggests that, as long as the constraint allows it, the minimizer of $K_n$ prefers to remain low-dimensional and to have coordinates grouped into equal pairs.  This motivates the following conjectural continuation beyond the four-dimensional endpoint.
	
	\begin{conjecture}\label{conj:six}
		Let $n\geq6$. For $d_{\max}(4)<\tau\leq d_{\max}(6)$, the minimizer of $K_n(\omega)$ on $\mathcal{C}_n(\tau)$ has the form
		\[
		\omega=(r,r,s,s,t,t,0,0,\ldots).
		\]
		More explicitly, the minimizer is conjectured to be
		\[
		(r,r,r,r,s,s,0,\ldots,0),
		\qquad
		d_{\max}(4)\leq\tau\leq\sqrt{2-\frac{46}{\sqrt{1107}}},
		\]
		and
		\[
		(r,r,t,t,t,t,0,\ldots,0),
		\qquad
		\sqrt{2-\frac{46}{\sqrt{1107}}}\leq\tau\leq d_{\max}(6),
		\]
		where
		\[
		r=\frac{2-\tau^2}{2\sqrt2},
		\qquad
		s=\sqrt{\frac12-2r^2},
		\qquad
		t=\frac12\sqrt{1-2r^2}.
		\]
	\end{conjecture}

\subsection{Common coefficient estimates}\label{sec:coefficient-estimates}
We use Haagerup's function
	\begin{equation}\label{eq:F-new}
	F(s)=\frac{2}{\sqrt{\pi s}}
	\frac{\Gamma((s+1)/2)}{\Gamma(s/2)},\qquad s>0.
	\end{equation}
	It is increasing, satisfies $F(4)=3/4$, and Haagerup's refined $L_1$-Khintchine estimate states that
	\begin{equation}\label{eq:haagerup-refined-new}
	K_n(\omega)\geq\sum_{i=1}^n\omega_i^2F(\omega_i^{-2})
	\end{equation}
	for every unit vector $\omega$ \cite{haagerup1981best}.  A summand corresponding to $\omega_i=0$ is interpreted by continuity.

	For the shared estimates used below, put
	\[
	C=\sqrt{\frac2\pi},\qquad
	r(u)=\sqrt{\frac{8+4u^2}{8+8u^2+3u^4}}\quad(u\geq0).
	\]

	\begin{lemma}[A sharper Wallis-ratio estimate]\label{lem:bks-wallis}
		For the function $F$ in \eqref{eq:F-new} and every $s>0$,
		\begin{equation}\label{eq:bks-wallis}
		F(s)\geq C\left(1+\frac1{2s}
		+\frac1{8s(s+1/2)}\right)^{-1/2}.
		\end{equation}
		Equivalently, $F(u^{-2})\geq Cr(u)$ for $u>0$.
	\end{lemma}

	The proof is given in Appendix~\ref{app:proof-lem-bks-wallis}.

	\begin{lemma}[Coefficient-dependent bounds]\label{lem:bks-coefficients}
		Let $\omega$ be a proper unit vector, and put $a=\omega_1$ and
		$b=\omega_2$, adding a zero coordinate if necessary. Then
		\begin{align}
		K_n(\omega)&\geq a,\label{eq:bks-largest}\\
		K_n(\omega)&\geq J(a,b):=\frac{a+b}{2}
		+\frac{\sqrt{1-2ab}}{2\sqrt2},\label{eq:bks-J}\\
		K_n(\omega)&\geq H(a,b):=
		C\bigl[a^2r(a)+(1-a^2)r(b)\bigr].\label{eq:bks-H}
		\end{align}
	\end{lemma}

	The proof is given in Appendix~\ref{app:proof-lem-bks-coefficients}.

	\begin{lemma}[A finite coefficient certificate]\label{lem:bks-certificate}
		Put $k=7423/10000$. If
		\[
		\frac12\leq a\leq k,\qquad 0\leq b\leq a,\qquad a+b\leq1,
		\]
		then $\max\{J(a,b),H(a,b)\}>k$.
	\end{lemma}

	The proof is given in Appendix~\ref{app:proof-lem-bks-certificate}.

\begin{corollary}[Scalar separation]\label{cor:scalar-separation}
Let $v$ be a unit vector, and let $v_{(1)}\geq v_{(2)}$ be its two
largest absolute coordinates, adjoining a zero coordinate if necessary.
If $v_{(1)}+v_{(2)}\leq1$, then
\[
 K(v)\geq k_0:=\frac{7423}{10000}.
\]
\end{corollary}
The proof is given in Appendix~\ref{app:scalar-separation}.

\subsection{The six-coordinate profile and the phase transition}
We use the following four-coordinate decomposition and localization estimate. The supporting proofs are collected in Appendix~\ref{app:supporting}.
\begin{definition} \label{def:kappa}
		For $c\in\mathbb{R}$ and $y=(y_1,\ldots,y_k)$, define
		\[
		\kappa(c;y)=\E\left|c\epsilon_0+\sum_{j=1}^k y_j\epsilon_j\right|,
		\]
		where $\epsilon_0,\epsilon_1,\ldots,\epsilon_k$ are independent Rademacher random variables.
	\end{definition}
	
	\begin{lemma}[Four-coordinate decomposition]\label{lem:four-decomposition}
		Let $y=(\omega_5,\ldots,\omega_n)$ and define
		\[
		q=\omega_1+\omega_2,
		\qquad
		u=\omega_3+\omega_4,
		\qquad
		d=\omega_1-\omega_2,
		\qquad
		e=\omega_3-\omega_4.
		\]
		Then
		\[
		K_n(\omega)
		\geq \frac38(q+u)
		+\frac18\bigl[\kappa(q-u;y)+\kappa(d+e;y)+\kappa(d-e;y)\bigr].
		\]
		Moreover, equality holds in a sufficiently small neighborhood of
		\[
		\omega_0=\left(\frac12,\frac12,\frac12,\frac12,0,0,\ldots\right),
		\]
		with the size of the neighborhood depending on $n$.
	\end{lemma}
	
	The proof is given in Appendix~\ref{app:proof-lem-four-decomposition}.

	For a proper vector, put
	\[
	x=q-u,
	\qquad
	t=\|y\|_2.
	\]
	Then $d,e\geq0$, and
	\[
	x-(d+e)=2(\omega_2-\omega_3)\geq0.
	\]
	If $\|\omega\|_2=1$, then
	\begin{equation}\label{eq:t-norm}
		t^2=1-q^2+qx-\frac{x^2+d^2+e^2}{2}.
	\end{equation}
	
	\begin{lemma}\label{lem:phi}
		Let $\omega\in\mathcal{C}_n(\tau)$ with $\tau> d_{\max}(4)$.  Then $q<1$.  Define
		\[
		A_q(x)=1-q^2+qx,
		\qquad
		B_q(x)=1-q^2+qx-x^2,
		\]
		and
		\[
		\Phi_q(x)=\frac34q-\frac38x
		+\frac{1}{8\sqrt2}\left(\sqrt{A_q(x)}+2\sqrt{B_q(x)}\right).
		\]
		Then
		\[
		K_n(\omega)\geq\Phi_q(x).
		\]
		Moreover, $\Phi_q'(0)>0$ if and only if $q>2\sqrt2/3$, and $\Phi_q''(x)\leq0$ on every interval where $A_q(x)$ and $B_q(x)$ are positive.
	\end{lemma}
	
	The proof is given in Appendix~\ref{app:proof-lem-phi}.

\begin{lemma}[Localization of competing vectors]\label{lem:six-localization}
Let $n\geq6$ and let $\omega\in\mathbb R^n$ be a proper unit vector with
$\sqrt{483}/22\leq q=\omega_1+\omega_2<1$, and put
\[
 A(q)=\frac{3q}{4}+\frac{3}{8\sqrt2}\sqrt{1-q^2},
 \qquad x=q-\omega_3-\omega_4.
\]
If $K_n(\omega)\leq A(q)$, then $0\leq x<1/12$.
\end{lemma}
The proof is given in Appendix~\ref{app:six-localization}.

\begin{theorem}[A dimension-free global six-coordinate interval]\label{th:six-local}
Let $n\geq6$ and set
\[
 d_0=d_{\max}(4)=\sqrt{2-\sqrt2},\qquad
 d_1=\sqrt{2-\frac{\sqrt{966}}{22}}.
\]
For every $\tau\in[d_0,d_1]$, the unique minimizer of $K_n$ on
$\mathcal C_n(\tau)$ is
\[
 \omega_\tau^{(6)}=(r,r,r,r,s,s,0,\ldots,0),\qquad
 r=\frac{2-\tau^2}{2\sqrt2},\quad
 s=\sqrt{\tau^2-\frac{\tau^4}{4}-\frac12}.
\]
In particular, writing $q=(2-\tau^2)/\sqrt2$, the minimum is
\[
 \min_{\omega\in\mathcal C_n(\tau)}K_n(\omega)
 =\frac{3q}{4}+\frac{3}{8\sqrt2}\sqrt{1-q^2}.
\]
The upper endpoint is $d_1=0.7663211842\ldots$, independent of $n$.
\end{theorem}

\begin{proof}
The endpoint $q=1$ follows from Theorem~\ref{thm:four}. For
$\sqrt{483}/22\leq q<1$, write $h=\sqrt{1-q^2}$ and let
$A(q)=3q/4+3h/(8\sqrt2)$. The displayed candidate is proper and
attains $A(q)$: its first four coefficients are $q/2$, its last two
nonzero coefficients are $h/\sqrt2$, and $q>\sqrt2h$.

It suffices to consider a competing vector with $K_n(\omega)\leq A(q)$.
By Lemma~\ref{lem:six-localization},
$x=q-\omega_3-\omega_4$ lies in $[0,1/12)$. Lemma~\ref{lem:phi}
gives $K_n(\omega)\geq\Phi_q(x)$. Here $q>99/100$ and
$h^2<1/100$. For $0<x\leq1/12$, both radicands in $\Phi_q$ are
positive and at most $1/100+1/12=7/75$, whence
\[
 \Phi_q'(x)=-\frac38+\frac1{8\sqrt2}
 \left(\frac{q}{2\sqrt{h^2+qx}}+
       \frac{q-2x}{\sqrt{h^2+qx-x^2}}\right)
 \geq-\frac38+\frac{791/600}{8\sqrt2\sqrt{7/75}}>0.
\]
The final strict inequality is certified by
$(791/600)^2-18(7/75)=20881/360000>0$.
Thus $K_n(\omega)\geq\Phi_q(x)\geq\Phi_q(0)=A(q)$, with equality
only if $x=0$. Equality then forces the first four coefficients to be
$q/2$, and equality in the tail Khintchine inequality forces precisely
two equal nonzero tail coefficients. Normalization gives
$\omega=\omega_\tau^{(6)}$, as claimed.
\end{proof}

Theorem~\ref{th:six-local} gives global minimality on a common interval
for all $n\geq6$. Together with Theorem~\ref{thm:four},
it gives the phase transition at $d_{\max}(4)$.
\begin{corollary}\label{cor:cusp}
For every $n\geq6$, as $\epsilon\downarrow0$,
\[
 \min_{\omega\in\mathcal C_n(d_{\max}(4)+\epsilon)}K_n(\omega)
 =\frac34+\frac{3\sqrt[4]{4-2\sqrt2}}8\sqrt\epsilon+O(\epsilon).
\]
\end{corollary}
The proof is given in Appendix~\ref{app:cusp}.

\section{An improved linear Khintchine inequality}\label{sec:robust}

We supplement the coefficient bounds with a sixth-moment estimate and
some scalar properties of $H$. Their proofs are given in
Appendix~\ref{app:robust-auxiliary}.

\begin{lemma}[Sixth-moment coefficient bound]\label{lem:robust-moment}
Let $\omega$ be a proper unit vector, $a=\omega_1>0$, and $b=\omega_2$.
Define
\begin{align}
 P_0(a)&=a+\frac{3(1-a^2)(-58a^4+63a^2-15)}{250a^5},
 \label{eq:robust-P0}\\
 P(a,b)&=P_0(a)+\frac{3b^2(1-a^2)}{250a^5}
       \min\{0,30-52a^2-16b^2\}.\label{eq:robust-P}
\end{align}
Then $K_n(\omega)\geq P(a,b)$, and $P_0$ is strictly increasing on
$(0,\infty)$.
\end{lemma}

\begin{lemma}[Curvature and endpoint estimates]\label{lem:robust-curvature}
For $0\leq u\leq1$, the function $r$ satisfies $r'(u)\leq0$ and
$r''(u)<0$; moreover, $r''(u)<-1/3$ for $0\leq u\leq1/2$.
For each $q\in[0,1]$, the function $a\mapsto H(a,q-a)$ is concave on
$[q/2,q]$. Put
\[
 \alpha=\frac{6817}{10000},\qquad
 T(q)=\frac1{\sqrt2}+\frac1{20}\sqrt{2-\sqrt2q}.
\]
The following strict inequalities hold, with all finite decimals
interpreted as exact rational numbers:
\begin{alignat*}{3}
 H(1/2,0)&>.7858,&\qquad H(1/2,1/2)&>.74993,\\
 H(\alpha,0)&>.75817,& H(\alpha,1-\alpha)&>.74750,\\
 T(0)&<.77783,& T(\alpha)&<.75801,\qquad T(1)<.74541,\\
 P_0(\alpha)&>.75811,& P_0(.74)&>.7915.
\end{alignat*}
\end{lemma}

The additional tail information in Lemma~\ref{lem:robust-moment}
allows us to prove the following dimension-free robust estimate.

\begin{theorem}[Improved linear Khintchine inequality]\label{th:explicit-robust}
For every $n\geq2$ and every unit vector $\omega\in\mathbb R^n$,
\begin{equation}\label{eq:robust-20}
 K_n(\omega)\geq\frac1{\sqrt2}+\frac1{20}\tau(\omega),
\end{equation}
where $\tau(\omega)$ is the distance to the extremal set defined in
\eqref{eq:proper}, after making $\omega$ proper by signs and permutations.
\end{theorem}

\begin{proof}
Assume $\omega$ proper and put $a=\omega_1$, $b=\omega_2$,
$q=a+b$, and $\tau=\sqrt{2-\sqrt2q}$.
If $q\geq1$, Lemma~\ref{lem:bks-coefficients} and $ab\leq q^2/4$ give
\[
 K_n(\omega)\geq\frac q2+\frac{\sqrt{2-q^2}}4
 =\frac1{\sqrt2}+\tau
   \frac{\sqrt{1-\tau^2/4}-\tau}{2\sqrt2}.
\]
The last coefficient decreases for $0\leq\tau\leq d_{\max}(4)$,
and at the right endpoint equals
\[
 \frac{3/4-1/\sqrt2}{\sqrt{2-\sqrt2}}>\frac1{20}.
\]
Thus it remains to prove $K_n(\omega)\geq T(q)$ when $q\leq1$.
In this range $\tau\geq\sqrt{2-\sqrt2}>3/4$, so
\begin{equation}\label{eq:robust-target-curvature}
 -T''(q)=\frac1{40\tau^3}<\frac8{135}<\frac1{16}.
\end{equation}

First suppose $a\leq\alpha=6817/10000$.
By the concavity in Lemma~\ref{lem:robust-curvature}, it is enough to
verify $H(a,q-a)>T(q)$ at the two endpoints of
$[q/2,\min\{q,\alpha\}]$.
At $a=q/2$, the difference $Cr(q/2)-T(q)$ is concave:
its first term has second derivative less than $-C/12<-1/16$,
where $C=\sqrt{2/\pi}>3/4$, and
\eqref{eq:robust-target-curvature} applies. Both endpoint values are
positive, since
\[
 C>T(0),\qquad Cr(1/2)>.74993>.74541>T(1).
\]
For the endpoint $a=q\leq\alpha$, set $H_0(q)=H(q,0)$.
This function decreases. When $q\leq1/2$, therefore,
\[
 H_0(q)\geq H_0(1/2)>.7858>.77783>T(0)\geq T(q).
\]
When $1/2\leq q\leq\alpha$, the bounds $r''<0$ and
$r'(1/2)<-895/3912<-1/5$ yield
\[
 H_0''(q)=C\bigl[2(r(q)-1)+4qr'(q)+q^2r''(q)\bigr]
 <-\frac3{10}.
\]
Consequently $H_0-T$ is concave, and its endpoint values are positive
by the preceding estimate and
$H_0(\alpha)>.75817>.75801>T(\alpha)$.
For the remaining endpoint $a=\alpha$, write $b=q-\alpha$.
Since $\alpha^2<1/2$ and $0\leq b\leq1-\alpha<1/2$,
\[
 H_{bb}(\alpha,b)=C(1-\alpha^2)r''(b)<-\frac18.
\]
Again $H(\alpha,b)-T(\alpha+b)$ is concave. It is positive at $b=0$
as just proved, and at $b=1-\alpha$ because
$H(\alpha,1-\alpha)>.74750>.74541>T(1)$.
This settles $a\leq\alpha$.

For $\alpha\leq a\leq.74$, use $b\leq1-a$ to obtain
\[
 30-52a^2-16b^2\geq14+32a-68a^2\geq.4432>0.
\]
Lemma~\ref{lem:robust-moment} gives
\[
 K_n(\omega)\geq P_0(a)\geq P_0(\alpha)
 >.75811>.75801>T(\alpha)\geq T(q).
\]
For $.74\leq a\leq.78$, we have $b\leq.26<.3$, $1-a^2<1/2$,
$30-52a^2-16b^2>-3$, and $250a^5>50$.
The correction in \eqref{eq:robust-P} is greater than $-.0081$;
hence
\[
 K_n(\omega)>P_0(.74)-.0081>.7834>.77783>T(0)\geq T(q).
\]
Finally, if $a\geq.78$, use $K_n(\omega)\geq a>T(0)\geq T(q)$.
\end{proof}

The four-coordinate unit vector $(1/2,1/2,1/2,1/2)$ gives the upper bound
\begin{equation}\label{eq:opt-const}
c_{\mathrm{opt}}\leq c_*:=\frac{3/4-1/\sqrt2}{\sqrt{2-\sqrt2}}
 =0.0560426911\ldots.
\end{equation}
Thus Theorem~\ref{th:explicit-robust} places the optimal universal slope in
$[1/20,c_*]$. The sharp four-coordinate profile proves the conjectured
slope for $\tau\leq d_{\max}(4)$; the remaining regime is open.
\begin{conjecture}\label{conj:optimal-slope}
For every unit vector $\omega$,
$K_n(\omega)\geq1/\sqrt2+c_*\tau(\omega)$, with $c_*$ as in
\eqref{eq:opt-const}.
\end{conjecture}

\section{Improved Bound on the Low-degree Fourier Weights of Boolean LTFs}\label{section:LTF}

Now we turn our attention to the low-degree Fourier weights of the Boolean Linear Threshold Functions (LTFs). Boolean functions are defined over the Hamming cube $\{-1, +1 \}^n$, and acts as important tools in computational complexity and learning theory. To start with, we recall some basics of analysis of Boolean functions.

Let $\mathcal{B}_n(\mathbb{R})=\{h:\{-1,1\}^n\to\mathbb{R}\}$ be the set of boolean functions where $n\in\mathbb{Z}_+$. Endow $\{-1,1\}^n$ with probablity measure $\mu_n$ such that $\mu_n(x)=\frac{1}{2^n}$ for every $x\in\{-1,1\}^n$. For every  $f,g\in\mathcal{B}_n(\mathbb{R})$, define its inner product $\langle f,g\rangle=\int fgd\mu_n$.
	
	For every $S\subseteq [n]:=\{1,2,\dots,n\}$, define $x_S:=\prod_{i\in S}x_i$. Then $\{x_S\}_{S\subseteq [n]}$ becomes an orthonormal basis for $\mathcal{B}_n(\mathbb{R})$. Therefore, for every $h\in\mathcal{B}_n(\mathbb{R})$, we can write $h=\sum_{S\subseteq[n]}\hat{h}(S)x_S$, where $\hat{h}(S)=\langle h, x_S\rangle$ is called the Fourier coefficient of $h$.
	
	We say $h\in\mathcal{B}_n(\mathbb{R})$ is a \textbf{linear threshold function} (LTF) if $h$ can be represented as
	$$h(x)=\sgn(\omega_0+\omega_1x_1+\omega_2x_2+\dots+\omega_nx_n)$$
	for some constants $\omega_0,\omega_1,\dots,\omega_n\in\mathbb{R}$, where $\sgn(x):=\begin{cases}
		1 & x\geq 0, \\
		-1 & x<0.
	\end{cases}$
	
	For every $h\in\mathcal{B}_n(\mathbb{R})$, we define 
	$$\W^{\leq 1}(h):=\sum_{\substack{S\subseteq [n] \\ |S|\leq 1}}\hat{h}(S)^2.$$
	We can see that $\W^{\leq 1}(h)=\|h^{\leq 1}\|_2^2$ with $h^{\leq 1}:=\sum_{|S|\leq 1}\hat{h}(S)x_S$.
	
	The following conjecture was first posed by Benjamini, Kalai, and Schramm \cite[Remark~3.7]{benjamini1999noise}; see also \cite[Conjecture~5.3]{odonnell2014analysis}.
	
	\begin{conjecture}[Benjamini--Kalai--Schramm]\label{W1-conj}
		Let $h:\{-1,1\}^n\to\{-1,1\}$ be an LTF. Then $\W^{\leq 1}(h)\geq\frac{2}{\pi}$.
	\end{conjecture}

    Beyond the classical bound $\mathbf W^{\leq1}[\mathrm{LTF}]\geq1/2$, De, Diakonikolas, and Servedio \cite{de2016robust} proved the existence of a universal $\epsilon>0$ such that $\mathbf W^{\leq1}[\mathrm{LTF}]\geq1/2+\epsilon$. We give the explicit bound in Theorem~\ref{th:explicit-ltf}.
	We first record reductions and a conditioned Chow estimate; their proofs are in Appendix~\ref{app:ltf-support}.
	
	\begin{observation}[Tie-free homogenization]\label{obs:homogenization}
For every Boolean LTF $f$ on a finite cube, there is a tie-free homogeneous
LTF $h$ on a cube with one additional coordinate such that
\begin{equation}\label{eq:homogenized-weight}
\W^{\leq1}[h]=\W^{\leq1}[f].
\end{equation}
Here tie-free means that the defining linear form is nonzero everywhere
on the cube.
\end{observation}
The proof is given in Appendix~\ref{app:homogenization}.

\begin{observation}[Symmetry reduction]\label{obs:symmetry}
Sign changes and coordinate permutations preserve the low-degree Fourier
weight and the tie-free property. A homogeneous representation may therefore
be chosen with a proper unit coefficient vector.
\end{observation}
The proof is given in Appendix~\ref{app:symmetry}.

\begin{observation}[Chow projection identity]\label{obs:chow-projection}
Let $h(x)=\sgn(\omega\cdot x)$ be a tie-free homogeneous LTF with
$\|\omega\|_2=1$, and let $\chi=(\widehat h(1),\ldots,\widehat h(n))$.
Then
\begin{equation}\label{eq:chow-defect-new}
\W^{\leq1}[h]=K_n(\omega)^2+
\|\chi-K_n(\omega)\omega\|_2^2\geq K_n(\omega)^2.
\end{equation}
\end{observation}
The proof is given in Appendix~\ref{app:chow-projection}.

For the following lemma, put $q=\omega_1+\omega_2$,
$p=\omega_1-\omega_2$, and $S=\sum_{i=3}^n\omega_i x_i$.
For $t\geq0$, define
\[
 b(t)=\E\sgn(S+t),\qquad
 c(t)=\E[(x_3,\ldots,x_n)\sgn(S+t)],\qquad
 A=b(q),\quad D=b(p).
\]
Symmetry and absence of ties give $A=\mathbb P(|S|<q)$ and
$D=\mathbb P(|S|<p)$.
\begin{lemma}\label{lem:close-to-two-coordinate}
		Let
		\[
		h(x)=\sgn\!\left(\sum_{i=1}^n\omega_i x_i\right),
		\qquad
		\|\omega\|_2=1,\qquad
		\omega_1\ge\cdots\ge\omega_n\ge0,
		\]
		and suppose that
		\[
		\sum_{i=1}^n\omega_i x_i\ne0
		\qquad\text{for every }x\in\{-1,1\}^n.
		\]
		With \(q,p,S,A,D\) defined above, one has
		\[
		W^{\le1}[h]\ge \frac{A^2}{2}+\frac18.
		\]
	\end{lemma}
	
	The proof is given in Appendix~\ref{app:proof-lem-close-to-two-coordinate}.

	\begin{theorem}[Explicit low-degree Fourier weight of LTFs]\label{th:explicit-ltf}
		Let
		\[
		\begin{aligned}
		\mathbf W^{\leq1}[\mathrm{LTF}]
		=\inf\bigl\{\mathbf W^{\leq1}[f]:{}&
		f\text{ is a Boolean linear threshold function}\\
		&\text{on a finite cube}\bigr\}.
		\end{aligned}
		\]
		Then
		\begin{equation}\label{eq:explicit-ltf-final}
		\mathbf W^{\leq1}[\mathrm{LTF}]
		\geq\left(\frac{7423}{10000}\right)^2
		=\frac{55100929}{100000000}=0.55100929.
		\end{equation}
	\end{theorem}

	\begin{proof}
		Constant LTFs have $\mathbf W^{\leq1}=1$, so let $f$ be nonconstant.
		By homogenization and \eqref{eq:homogenized-weight}, it suffices to
		consider a tie-free homogeneous LTF $h$ with a proper unit
		coefficient vector $\omega$. Put
		\[
		a=\omega_1,\qquad b=\omega_2,\qquad q=a+b,\qquad
		k=\frac{7423}{10000}.
		\]

		First suppose that $q\geq4/3$. Since
		$\omega_1^2+\omega_2^2\geq q^2/2$, Chebyshev's inequality gives
		\[
		A=\mathbb P(|S|<q)
		\geq1-\frac{\E S^2}{q^2}
		\geq1-\frac{1-q^2/2}{q^2}
		=\frac32-\frac1{q^2}\geq\frac{15}{16}.
		\]
		Lemma~\ref{lem:close-to-two-coordinate} therefore gives
		\[
		\W^{\leq1}[h]\geq\frac12\left(\frac{15}{16}\right)^2+\frac18
		=\frac{289}{512}>k^2.
		\]

		Next suppose that $1\leq q\leq4/3$. Since $ab\leq q^2/4$,
		\eqref{eq:bks-J} gives
		\[
		K_n(\omega)\geq\frac q2+\frac{\sqrt{2-q^2}}4.
		\]
		The right-hand side is concave on $[1,4/3]$. Its endpoint values
		are $3/4$ and $2/3+\sqrt2/12>3/4$, respectively. Hence
		$K_n(\omega)\geq3/4>k$, and \eqref{eq:chow-defect-new} applies.

		Finally suppose that $q\leq1$. Corollary~\ref{cor:scalar-separation}
gives $K_n(\omega)\geq k$.

All cases yield $\W^{\leq1}[h]\geq k^2$ by
		\eqref{eq:chow-defect-new}. Homogenization transfers the bound
		to every Boolean LTF, proving \eqref{eq:explicit-ltf-final}.
	\end{proof}

\section{Conclusion}
We determine the sharp four-coordinate profile and an explicit interval, uniform in dimension, on which the six-coordinate family is the unique global minimizer. These results establish a square-root transition at $d_{\max}(4)$. The global coefficient $1/20$ in the linear Khintchine inequality and the Fourier-weight
bound $0.55100929$ follow from coefficient-dependent estimates, a
sixth-moment tail inequality, and conditioned Chow parameters.
The optimal robust slope and the $2/\pi$ Fourier-weight conjecture remain open. 
Analyzing properties of minimizing profiles in high dimensions is a promising direction for future research. 
%

\section*{Statement on AI Usage, Human Effort and Revision History}

    The core results in this paper are AI-assisted rather than directly produced by AI. The first version of this work was made public before the release of GPT-6, and some of the results in the paper were completed before the release of GPT-5.6. In the first version, AI was used primarily for background knowledge, numerical simulations and idea verification. For the first version, the AI systems used are ChatGPT and DeepSeek, supported by the harness platforms Cursor, Codex, and SeedProver, as well as an open-source harness system Conjecta\footnote{https://github.com/conjecta/conjecta}, written by J. Yu, who was then a graduate student at Fudan University. 
    
    Later the authors used GPT-6 to refine and optimize some of the arguments in the second version of this article. In particular, GPT-6 gives the current version of Theorem 2, which is an improvement over the phase transition result in the first version. Also, GPT-6 gives the current bound on the BKS constant, which is an improvement over the bound in the first version.

    
    The authors take full responsibility of the content in the paper. 

    \section*{Acknowledgement} 
    
    Part of the work was done while Xuan was an intern at ByteDance. 

\clearpage
\appendix
\section{Proofs of supporting results}\label{app:supporting}
The proofs omitted from the main text are collected here. Proofs for the additional properties stated in Appendix~\ref{app:additional} are included at the end of this appendix.

\subsection{Proof of Lemma~\ref{lem:bks-wallis}}\label{app:proof-lem-bks-wallis}
\begin{proof}[Proof of Lemma~\ref{lem:bks-wallis}]
		For $x>0$, define
		\[
		G(x)=\frac{\Gamma(x+1/2)}{\Gamma(x)},\qquad
		D(x)=x+\frac14+\frac1{32x+8},\qquad
		L(x)=\frac{x}{\sqrt{D(x)}}.
		\]
		Direct expansion gives
		\[
		(x+1)^2D(x)-(x+\tfrac12)^2D(x+1)
		=\frac9{32(4x+1)(4x+5)}>0.
		\]
		Since $G(x+1)/G(x)=(x+1/2)/x$, it follows that
		\[
		\frac{G(x+1)/L(x+1)}{G(x)/L(x)}
		=\frac{x+1/2}{x+1}\sqrt{\frac{D(x+1)}{D(x)}}<1.
		\]
		For fixed $x>0$, Stirling's formula gives
		$G(x+m)/L(x+m)\to1$ as $m\to\infty$. This sequence is decreasing,
		so $G(x)>L(x)$. Substituting $x=s/2$ into \eqref{eq:F-new} gives
		\[
		F(s)>C\sqrt{\frac{x}{D(x)}}
		=C\left(1+\frac1{2s}+\frac1{8s(s+1/2)}\right)^{-1/2}.
		\]
		The equivalent form follows by setting $s=u^{-2}$.
	\end{proof}

\subsection{Proof of Lemma~\ref{lem:bks-coefficients}}\label{app:proof-lem-bks-coefficients}
\begin{proof}[Proof of Lemma~\ref{lem:bks-coefficients}]
		Symmetry and Jensen's inequality give
		$K_n(\omega)=\E|a+\sum_{i\geq2}\omega_i x_i|\geq a$.
		Conditioning on $x_1,x_2$, and using the symmetry of $S$, gives
		\[
		K_n(\omega)=\tfrac12\E|q+S|+\tfrac12\E|p+S|.
		\]
		The first expectation is at least $q=a+b$. For an independent
		Rademacher variable $z$, the sharp Khintchine inequality gives
		\[
		\E|p+S|=\E|pz+S|
		\geq\frac1{\sqrt2}\sqrt{p^2+\E S^2}
		=\frac1{\sqrt2}\sqrt{1-2ab},
		\]
		which proves \eqref{eq:bks-J}.
		Finally, $r$ is decreasing, since, with $t=u^2$,
		\[
		\frac{d}{dt}\frac{8+4t}{8+8t+3t^2}
		=-\frac{32+48t+12t^2}{(8+8t+3t^2)^2}<0.
		\]
		Apply \eqref{eq:haagerup-refined-new} and
		Lemma~\ref{lem:bks-wallis}, then use $\omega_i\leq b$ for $i\geq2$
		and $\sum_{i\geq2}\omega_i^2=1-a^2$ to obtain \eqref{eq:bks-H}.
	\end{proof}

\subsection{Proof of Lemma~\ref{lem:bks-certificate}}\label{app:proof-lem-bks-certificate}
\begin{proof}[Proof of Lemma~\ref{lem:bks-certificate}]
		On $a,b\geq0$ and $a+b\leq1$, $J$ is increasing in each coordinate.
		Indeed,
		\[
		\partial_aJ=\frac12-\frac{b}{2\sqrt2\sqrt{1-2ab}}>0,
		\]
		because $4ab+b^2\leq4b-3b^2\leq4/3<2$; the other derivative
		follows by symmetry. On $0\leq b\leq a\leq1$, $H$ is decreasing
		in each coordinate, since
		\[
		\partial_aH=C\{2a[r(a)-r(b)]+a^2r'(a)\}\leq0,
		\qquad
		\partial_bH=C(1-a^2)r'(b)\leq0.
		\]

		Set
		\[
		c_-=\frac{7978845}{10^7},\qquad
		d_-=\frac{3535533}{10^7}.
		\]
		The classical enclosure $\pi<355/113$ and the exact rational
		inequalities
		\[
		2-\frac{355}{113}c_-^2
		=\frac{12203369}{90400000000000}>0,
		\qquad
		1-8d_-^2=\frac{6405911}{12500000000000}>0
		\]
		show that $c_-<C$ and $d_-<1/(2\sqrt2)$.
		For each row $(\ell,u,t,\alpha,\beta,\zeta)$ of
		Table~\ref{tab:bks-radicals}, direct multiplication of positive
		denominators verifies
		\[
		\alpha^2<\frac{8+4u^2}{8+8u^2+3u^4},\qquad
		\beta^2<\frac{8+4t^2}{8+8t^2+3t^4},\qquad
		\zeta^2<1-2\ell t.
		\]
		All finite decimals in the tables denote exact rational numbers.
		Consequently
		\[
		J(\ell,t)>\frac{\ell+t}{2}+d_-\zeta,
		\qquad
		H(u,t)>c_-\bigl[u^2\alpha+(1-u^2)\beta\bigr].
		\]
		Table~\ref{tab:bks-endpoints} records rational lower bounds for
		these expressions, rounded down to seven decimal places.
		Every entry is greater than $k$.

		\begin{table}[htbp]
		\centering
		\small
		\setlength{\tabcolsep}{6pt}
		\begin{tabular}{cccccc}
		\hline
		$\ell$ & $u$ & $t$ & $\alpha<r(u)$ & $\beta<r(t)$ &
		$\zeta<\sqrt{1-2\ell t}$\\
		\hline
		0.5   & 0.63   & 0.47  & 0.907263 & 0.946627 & 0.728010\\
		0.63  & 0.71   & 0.294 & 0.884775 & 0.978647 & 0.793448\\
		0.71  & 0.738  & 0.144 & 0.876564 & 0.994829 & 0.891919\\
		0.738 & 0.7423 & 0.086 & 0.875289 & 0.998152 & 0.934378\\
		\hline
		\end{tabular}
		\caption{Exact rational lower bounds for the three radicals.}
		\label{tab:bks-radicals}
		\end{table}

		\begin{table}[htbp]
		\centering
		\small
		\begin{tabular}{cccc}
		\hline
		$[\ell,u]$ & $t$ & Lower bound for $J(\ell,t)$ &
		Lower bound for $H(u,t)$\\
		\hline
		$[0.5,0.63]$     & 0.47  & 0.7423903 & 0.7428332\\
		$[0.63,0.71]$    & 0.294 & 0.7425261 & 0.7430906\\
		$[0.71,0.738]$   & 0.144 & 0.7423409 & 0.7423650\\
		$[0.738,0.7423]$ & 0.086 & 0.7423524 & 0.7423943\\
		\hline
		\end{tabular}
		\caption{Certified endpoint bounds, all exceeding $k=0.7423$.}
		\label{tab:bks-endpoints}
		\end{table}

		The four intervals $[\ell,u]$ cover $[1/2,k]$. For a point in
		one such interval, if $b\geq t$, monotonicity gives
		$J(a,b)\geq J(\ell,t)>k$; this comparison stays in $a+b\leq1$.
		If $b\leq t$, then $\ell\geq t$ in each row, and monotonicity gives
		$H(a,b)\geq H(u,t)>k$. This proves the claim.
	\end{proof}

\subsection{Proof of Proposition~\ref{prop:dmax}}\label{app:proof-prop-dmax}
\begin{proof}[Proof of Proposition~\ref{prop:dmax}]
		For a proper unit vector, put $q=\omega_1+\omega_2$.  Since
		\[
		\tau^2=\|\omega-\omega^*\|_2^2=2-\sqrt2\,q,
		\]
		maximizing $\tau$ is the same as minimizing $q$.
		
		Fix $q$ and write $y=\omega_2$.  Then $0\leq y\leq q/2$, $\omega_1=q-y$, and all coordinates $\omega_3,\ldots,\omega_n$ are at most $y$.  Hence
		\[
		1=\|\omega\|_2^2
		\leq (q-y)^2+(n-1)y^2.
		\]
		Maximizing the right-hand side over $0\leq y\leq q/2$ gives
		\[
		1\leq q^2\max\left\{1,\frac n4\right\}.
		\]
		Thus $q\geq1$ for $2\leq n\leq4$, and $q\geq 2/\sqrt n$ for $n\geq4$.  These lower bounds are attained by $(1,0,\ldots,0)$ when $n\leq4$ and by the constant vector $(1/\sqrt n,\ldots,1/\sqrt n)$ when $n\geq4$.  Substituting the minimum value of $q$ into $\tau^2=2-\sqrt2 q$ gives the stated formula.
	\end{proof}

\subsection{Proof of Lemma~\ref{lem:four-decomposition}}\label{app:proof-lem-four-decomposition}
\begin{proof}[Proof of Lemma~\ref{lem:four-decomposition}]
		Condition on the last $n-4$ variables and average over the first four signs.  Among the $16$ sign patterns of $(\epsilon_1,\ldots,\epsilon_4)$, there are $10$ unbalanced patterns, namely those with a strict majority of one sign.  For such a pattern let
		\[
		A=\sum_{i=1}^4\omega_i\epsilon_i,
		\qquad
		Y=\sum_{i=5}^n\omega_i\epsilon_i.
		\]
		Since $\E Y=0$, for each fixed sign pattern $(\epsilon_1,\epsilon_2,\epsilon_3,\epsilon_4)$, Jensen's inequality gives
		\[
		\E_Y|A+Y|\geq |A|.
		\]
		Summing $|A|$ over the $10$ unbalanced patterns, the triangle inequality gives a lower bound $6(q+u)$, and therefore their total contribution to $K_n(\omega)$ is at least $3(q+u)/8$.
		
		The remaining $6$ sign patterns are balanced.  Their first-four-coordinate sums are
		\[
		\pm(q-u),\qquad \pm(d+e),\qquad \pm(d-e).
		\]
		Averaging each opposite pair gives respectively
		\[
		\frac18\kappa(q-u;y),
		\qquad
		\frac18\kappa(d+e;y),
		\qquad
		\frac18\kappa(d-e;y).
		\]
		This proves the global lower bound.
		
		If $\omega$ is sufficiently close to $\omega_0$, then for every unbalanced pattern the quantity $|A|$ is strictly larger than $\|y\|_1$.  In that case the sign of $A+Y$ cannot change as $Y$ varies, and hence $\E|A+Y|=|A|$.  Therefore the preceding lower bound is an equality in such a neighborhood.
	\end{proof}

\subsection{Proof of Lemma~\ref{lem:phi}}\label{app:proof-lem-phi}
\begin{proof}[Proof of Lemma~\ref{lem:phi}]
		Since $\omega\in\mathcal{C}_n(\tau)$,
		\[
		\tau^2=2-\sqrt2 q.
		\]
		Thus $\tau> d_{\max}(4)=\sqrt{2-\sqrt2}$ implies $q<1$.
		By the classical Khintchine inequality,
		\[
		\kappa(c;y)\geq\frac1{\sqrt2}\sqrt{c^2+t^2}.
		\]
		Using Lemma \ref{lem:four-decomposition} and the relation $u=q-x$, we obtain
		
		\begin{equation}\label{eq:8}
			K_n(\omega)\geq \frac34q-\frac38x+\frac1{8\sqrt2}\left(\sqrt{x^2+t^2}+\sqrt{(d+e)^2+t^2}+\sqrt{(d-e)^2+t^2}\right).
		\end{equation}
		
		Because $d,e\geq0$ and $d+e\leq x$, we have
		\[
		d^2+e^2\leq(d+e)^2\leq x^2.
		\]
		Together with \eqref{eq:t-norm}, this gives
		\[
		x^2+t^2=A_q(x)+\frac{x^2-d^2-e^2}{2}\geq A_q(x)
		\]
		and
		\[
		t^2=B_q(x)+\frac{x^2-d^2-e^2}{2}\geq B_q(x).
		\]
		Substituting these inequalities into the previous lower bound yields
		\[
		K_n(\omega)\geq\Phi_q(x).
		\]
		
		Finally,
		\[
		\Phi_q'(0)
		=-\frac38+\frac{3q}{16\sqrt2\sqrt{1-q^2}},
		\]
		which is positive exactly when $q>2\sqrt2/3$.  The inequality $\Phi_q''(x)\leq0$ follows from the concavity of the square-root function: $A_q$ is affine in $x$, while $B_q$ is concave in $x$ on its positivity interval.
	\end{proof}

\subsection{Scalar separation}\label{app:scalar-separation}
\begin{proof}[Proof of Corollary~\ref{cor:scalar-separation}]
By symmetry we may make $v$ proper. Write $a=v_1$, $b=v_2$ and
$k_0=7423/10000$. If $a\geq k_0$, use $K(v)\geq a$ from
Lemma~\ref{lem:bks-coefficients}. If $a\leq1/2$, Haagerup's refined
inequality~\eqref{eq:haagerup-refined-new}, the monotonicity of $F$ and
$F(4)=3/4$ give $K(v)\geq3/4>k_0$. Otherwise,
Lemma~\ref{lem:bks-certificate} gives
$K(v)\geq\max\{J(a,b),H(a,b)\}>k_0$.
\end{proof}

\subsection{Localization of six-coordinate competitors}\label{app:six-localization}
\begin{proof}[Proof of Lemma~\ref{lem:six-localization}]
Put $p=\omega_1-\omega_2$, $h=\sqrt{1-q^2}$,
$S=\sum_{i\geq3}\omega_i\epsilon_i$, and
\[
 v=(p,\omega_3,\ldots,\omega_n),\qquad
 T=\|v\|_2,\qquad T^2=1-\frac{q^2}{2}+\frac{p^2}{2}.
\]
Then $1/\sqrt2\leq T\leq1$, $q\leq\sqrt2T$ and $h\leq1/22$.
Conditioning on the first two signs gives
\[
 K_n(\omega)=\tfrac12\E|q+S|+\tfrac12\E|p+S|
 \geq\frac q2+\frac{K(v)}2.
\]
Our assumption $K_n(\omega)\leq A(q)$ therefore yields
\begin{equation}\label{eq:six-aux-upper}
 \frac{K(v)}T\leq\frac q{2T}+\frac{3h}{4\sqrt2T}
 \leq\frac1{\sqrt2}+\frac34h
 \leq\frac1{\sqrt2}+\frac3{88}<k_0.
\end{equation}
For the last strict comparison use $1/\sqrt2<177/250$ and the exact identity
\[
 k_0-\left(\frac{177}{250}+\frac3{88}\right)=\frac{23}{110000}>0.
\]
Let $\rho=\tau(v/T)$. Corollary~\ref{cor:scalar-separation} shows that
the two largest coordinates of $v/T$ sum to more than $1$, hence
$\rho<d_0:=d_{\max}(4)$. By Theorem~\ref{thm:four} and
\eqref{eq:six-aux-upper},
\[
 f(\rho):=\frac{\rho}{2\sqrt2}
 \left(\sqrt{1-\rho^2/4}-\rho\right)
 \leq\frac34h\leq\frac3{88}<\frac1{25}.
\]
The function $f$ is concave on $[0,d_0]$. Indeed, if
$R=1-z^2/4$, then
\[
 f''(z)=-\frac1{2\sqrt2}\left(2+\frac{z(6-z^2)}{8R^{3/2}}\right)<0.
\]
Its endpoint values on $[1/5,d_0]$ both exceed $1/25$:
$f(d_0)=3/4-1/\sqrt2>1/25$, and, for $0<z\leq1/5$,
\[
 f(z)\geq\frac{z}{2\sqrt2}
 \left(\sqrt{99/100}-1/5\right)>\frac z4.
\]
The latter strict inequality follows from
$\sqrt{99/100}>99/100$ and $1/\sqrt2<177/250$; in particular
$f(1/5)>1/20$. Thus $\rho<1/5$, and the same estimate gives
$\rho\leq3h$ (with $\rho=0$ allowed).

The nearest canonical pair cannot contain the $p$ coordinate. Every
tail coordinate $y$ is at most $\omega_2=(q-p)/2$, so
\[
 \frac{p+y}{T}\leq\frac{q+p}{2T}\leq1,
 \qquad
 4T^2-(q+p)^2=4-4q^2+(p-q)^2\geq0.
\]
A pair containing $p$ would therefore have squared distance at least
$2-\sqrt2$, contrary to $\rho<d_0$. Hence a closest pair is given by
$\omega_3,\omega_4$, and the distance formula gives
\[
 \omega_3+\omega_4=\frac{T}{\sqrt2}(2-\rho^2).
\]
Since $T\leq1$ and $q\leq\sqrt2T$, we conclude that
\[
 0\leq x=q-\sqrt2T+\frac{T\rho^2}{\sqrt2}
 \leq\frac{\rho^2}{\sqrt2}
 \leq\frac{9h^2}{\sqrt2}
 \leq\frac9{484\sqrt2}<\frac1{12}.
\]
\end{proof}

\subsection{The square-root expansion}\label{app:cusp}
\begin{proof}[Proof of Corollary~\ref{cor:cusp}]
		For $\omega_\tau^{(6)}$,
		\[
		K_n(\omega_\tau^{(6)})=\Phi_q(0)
		=\frac34q+\frac{3}{8\sqrt2}\sqrt{1-q^2},
		\qquad
		q=\frac{2-\tau^2}{\sqrt2}.
		\]
		Set $d_0=d_{\max}(4)=\sqrt{2-\sqrt2}$ and $\tau=d_0+\epsilon$.  Since $q=1-\sqrt2d_0\epsilon+O(\epsilon^2)$, we have
		\[
		\sqrt{1-q^2}=\sqrt{2\sqrt2d_0}\,\sqrt\epsilon+O(\epsilon^{3/2}).
		\]
		As $\sqrt{\sqrt2d_0}=\sqrt[4]{4-2\sqrt2}$, the desired expansion follows.
	\end{proof}
\subsection{Auxiliary estimates for the robust coefficient}
\label{app:robust-auxiliary}

\begin{proof}[Proof of Lemma~\ref{lem:robust-moment}]
We first prove the polynomial minorant
\begin{equation}\label{eq:robust-minorant}
 (x-1)_+\geq\frac3{25}(x^2-1)x^2(1-x^2/10)\qquad(x\geq0).
\end{equation}
For $0\leq x\leq1$ or $x\geq\sqrt{10}$ the right side is nonpositive.
In the remaining interval, division by $x-1$ reduces the claim to
$h(x)=250-30x^2-30x^3+3x^4+3x^5\geq0$.
With $y=x-12/5\geq-7/5$, exact expansion gives
\[
 h=\frac{2786}{3125}+\frac{144}{125}y+\frac{1362}{5}y^2
   +\frac{858}{5}y^3+39y^4+3y^5.
\]
For $y\geq0$ all terms are nonnegative. For $-7/5\leq y\leq0$,
the bracket multiplying $y^2$ is at least
\[
 \frac{1362}{5}-\frac{858}{5}\frac75-3\left(\frac75\right)^3
 =\frac{2991}{125}>20.
\]
Also $2786/3125>4/5$ and $(144/125)y\geq(6/5)y$.
Thus $h\geq4/5+(6/5)y+20y^2>0$, proving the minorant.

Write $S=\sum_{i\geq2}\omega_i\epsilon_i$, $s^2=1-a^2$,
$Q_4=\sum_{i\geq2}\omega_i^4$, and $Q_6=\sum_{i\geq2}\omega_i^6$.
Conditioning on $\epsilon_1$ gives $K_n(\omega)=\E\max\{a,|S|\}$.
Applying \eqref{eq:robust-minorant} to $|S|/a$ yields
\[
 K_n(\omega)\geq a+\frac3{250a^5}
 \left(11a^2\E S^4-10a^4\E S^2-\E S^6\right).
\]
The exact Rademacher moment identities
\[
 \E S^2=s^2,\qquad \E S^4=3s^4-2Q_4,\qquad
 \E S^6=15s^6-30s^2Q_4+16Q_6
\]
transform the bracket into
\[
 (1-a^2)(-58a^4+63a^2-15)+(30-52a^2)Q_4-16Q_6.
\]
Since $Q_6\leq b^2Q_4$ and $0\leq Q_4\leq b^2(1-a^2)$,
this proves \eqref{eq:robust-P}.
The derivative of $P_0$, multiplied by $250a^6$, is
$424u^3+363u^2-702u+225$, where $u=a^2$.
Putting $v=u-1/2\geq-1/2$, this is
\[
 424v^3+999v^2-21v+\frac{71}{4}
 \geq787v^2-21v+\frac{71}{4}>0.
\]
The last quadratic has negative discriminant, so $P_0$ increases.
\end{proof}

\begin{proof}[Proof of Lemma~\ref{lem:robust-curvature}]
Put $t=u^2$ and $Q(t)=(8+4t)(8+8t+3t^2)$.
Direct differentiation gives
\[
 r'(u)=-r(u)u\frac{32+48t+12t^2}{Q(t)},\qquad
 r''(u)=r(u)\frac{S(t)}{Q(t)^2},
\]
where $S(t)=-2048-5120t+4224t^3+2208t^4+288t^5$.
For $0\leq t\leq1$,
$S(t)\leq-2048+1600t\leq-448$.
For $0\leq t\leq1/4$, we have
\[
 r(u)\geq\frac{15}{16},\qquad -S(t)\geq2048+4800t,
 \qquad Q(t)\leq64+\frac{443}{4}t.
\]
Indeed $r(1/2)^2=144/163>(15/16)^2$; the bounds for $S$ and $Q$
follow by using $t^j\leq(1/4)^{j-1}t$.
Similarly $Q(t)^2\leq4096+(1103513/64)t$, whence
\[
 3r(u)(-S(t))-Q(t)^2
 \geq1664-\frac{239513}{64}t
 \geq\frac{186471}{256}>0.
\]
This proves $r''<-1/3$ on $[0,1/2]$. Also
$r'(1/2)=-(358/1467)r(1/2)<-895/3912$, using $r(1/2)>15/16$. The same derivative formula and
$r(1/2)>15/16$ give
$ r'(1/2)<-(15/16)(358/1467)=-895/3912<-1/5$.
For $b=q-a\leq a\leq1$,
\begin{align*}
 \frac1C\frac{d^2}{da^2}H(a,q-a)
 ={}&2[r(a)-r(b)]+4a[r'(a)+r'(b)]\\
 &+a^2r''(a)+(1-a^2)r''(b)\leq0,
\end{align*}
proving the required concavity.

For completeness, all endpoint comparisons admit the following exact
rational certificates. Put $C_-=.7978845$.
The classical bound $\pi<355/113$ and $C_-^2<226/355$ show $C_-<C$;
also $1/\sqrt2<.70711$ and $\sqrt2>1.41421$ by squaring.
The radical lower witnesses
\[
 r(.5)>.939912,\qquad r(.6817)>.892906,\qquad r(.3183)>.975028
\]
are verified by squaring and substituting the rational expression for
$r(u)^2$. Substitution into
$C_-[a^2r(a)+(1-a^2)r(b)]$ gives, respectively,
\[
\begin{array}{c|c}
 (a,b)&\text{lower bound for }H(a,b)\\ \hline
 (.5,0)&.7858\\
 (.5,.5)&.74993\\
 (.6817,0)&.75817\\
 (.6817,.3183)&.74750
\end{array}
\]
For the target, squaring verifies
\[
 \sqrt2<1.41422,\qquad
 \sqrt{2-\sqrt2\alpha}<1.0179,\qquad
 \sqrt{2-\sqrt2}<.766,
\]
where the lower bound $\sqrt2>1.41421$ is used for the last two.
Together with $1/\sqrt2<.70711$, these give the three asserted bounds
for $T$. Finally the rational substitutions in \eqref{eq:robust-P0}
give $P_0(.6817)>.75811$ and $P_0(.74)>.7915$.
All comparisons are strict and involve rational arithmetic after the
indicated squarings.
\end{proof}

\subsection{Supporting proofs for LTFs}\label{app:ltf-support}
\subsubsection{Tie-free homogenization}\label{app:homogenization}
\begin{proof}[Proof of Observation~\ref{obs:homogenization}]
Write $f(x)=\sgn(w_0+\sum_{i=1}^m w_ix_i)$ with $\sgn(0)=1$.
Since the cube is finite, increase $w_0$ by a sufficiently small positive
amount to obtain $w'_0$ without changing $f$ and so that
$w'_0+\sum_iw_ix_i$ never vanishes. Put
$h(y,x)=\sgn(w'_0y+\sum_iw_ix_i)$.
This form is tie-free: a zero for $y=-1$ would give a zero for $y=1$
after replacing $x$ by $-x$.
Moreover $h(1,x)=f(x)$ and $h(-1,x)=-f(-x)$. Averaging over $y$ gives
\[
\widehat h(y)=\widehat f(\emptyset),\qquad
\widehat h(i)=\widehat f(i)\ (1\leq i\leq m),\qquad
\widehat h(\emptyset)=0.
\]
The identity for the Fourier weights follows.
\end{proof}

\subsubsection{Symmetry reduction}\label{app:symmetry}
\begin{proof}[Proof of Observation~\ref{obs:symmetry}]
A coordinate sign change changes the sign of the corresponding first-level
Fourier coefficient, and a coordinate permutation permutes those coefficients.
Neither operation changes their squared sum or the constant coefficient.
They also preserve the set of values of the defining linear form.
Finally positive rescaling preserves its sign. These operations yield a
proper unit coefficient vector without affecting tie-freeness or the weight.
\end{proof}

\subsubsection{Chow projection}\label{app:chow-projection}
\begin{proof}[Proof of Observation~\ref{obs:chow-projection}]
Homogeneity and absence of ties make $h$ odd, hence balanced. Also
\[
\langle\chi,\omega\rangle=\E[h(x)(\omega\cdot x)]
=\E|\omega\cdot x|=K_n(\omega).
\]
The vectors $K_n(\omega)\omega$ and $\chi-K_n(\omega)\omega$ are
orthogonal. Pythagoras gives \eqref{eq:chow-defect-new}.
\end{proof}

\subsection{Proof of Lemma~\ref{lem:close-to-two-coordinate}}\label{app:proof-lem-close-to-two-coordinate}
\begin{proof}[Proof of Lemma~\ref{lem:close-to-two-coordinate}]
		Averaging over the four choices of $(x_1,x_2)$ gives
		\begin{equation}\label{eq:conditioned-chow-new}
			\W^{\leq 1}[h]
			=\frac{A^2+D^2}{2}
			+\frac14\|c(q)+c(p)\|_2^2.
		\end{equation}
		Indeed, the first two Chow parameters are 

		$$\hat{h}(\{1\})=\mathbb{E}\sgn\left(\omega_1+\omega_2x_2+S\right)=\mathbb{P}\left(\left|\omega_2x_2+S\right|\leq\omega_1\right)=\frac{A+D}{2},$$
		$$\hat{h}(\{2\})=\mathbb{E}\sgn\left(\omega_1x_1+\omega_2+S\right)=\mathbb{P}\left(\left|\omega_1x_1+S\right|\leq\omega_2\right)=\frac{A-D}{2},$$
		for $i\geq 3$,
		$$\hat{h}(\{i\})=\mathbb{E}\left[x_i\sgn(\omega_1x_1+\omega_2x_2+S)\right]=\frac{1}{2}\left(\mathbb{E}\left[x_i\sgn(q+S)\right]+\mathbb{E}\left[x_i\sgn(p+S)\right]\right),$$
		so the remaining Chow vector is $(c(q)+c(p))/2$, and we get \eqref{eq:conditioned-chow-new}.
		
		Each coordinate of $c(t)$ is nonnegative: for each $i\geq 3$, $\mathbb{E}[x_i\sgn(S+t)]=\mathbb{E}[\sgn(\omega_i+\sum_{\substack{j\geq3 \\ j\neq i}}\omega_jx_j+tx_i)]=\mathbb{P}(\|\sum_{\substack{j\geq3 \\ j\neq i}}\omega_jx_j+tx_i\|\leq\omega_i)\geq 0$. Hence
		\[
		\langle c(q),c(p)\rangle\geq0.
		\]
		The balanced homogeneous LTF
		\[
		(x_0,x_3,\ldots,x_n)
		\longmapsto \sgn(px_0+S)
		\]
		has first-level Chow vector $(D,c(p))$.  Applying the sharp Khintchine inequality and then Cauchy--Schwarz gives
		\[
		\frac1{\sqrt2}\sqrt{p^2+\sum_{i=3}^n\omega_i^2}
		\leq \E|px_0+S|
		\leq
		\sqrt{D^2+\|c(p)\|_2^2}
		\sqrt{p^2+\sum_{i=3}^n\omega_i^2}.
		\]
		Consequently,
		\[
		D^2+\|c(p)\|_2^2\geq\frac12.
		\]
		Using this and the coordinatewise nonnegativity in \eqref{eq:conditioned-chow-new}, we obtain the useful self-contained estimate
		\begin{equation}\label{eq:near-chow-bound}
			\W^{\leq 1}[h]
			\geq\frac{A^2}{2}+\frac18.
		\end{equation}
	\end{proof}

\subsection{Proofs for the additional properties}\label{app:additional-proofs}
In this subsection, $a(\omega)$ and the convention $\sgn(0)=0$ are as in Appendix~\ref{app:additional}.

\subsubsection{Proof of Proposition~\ref{prop:directional}}\label{app:proof-prop-directional}
\begin{proof}[Proof of Proposition~\ref{prop:directional}]
		For each fixed sign vector $\epsilon\in\{-1,1\}^n$, the one-sided derivative of
		\[
		\delta\mapsto |(\omega+\delta h)\cdot\epsilon|
		\]
		at $\delta=0$ is
		\[
		\sgn(\omega\cdot\epsilon)h\cdot\epsilon
		\quad\text{if }\omega\cdot\epsilon\neq0,
		\]
		and is $|h\cdot\epsilon|$ if $\omega\cdot\epsilon=0$.  Averaging over all sign vectors gives
		\begin{align*}
			DK_n(\omega;h)
			&=\E\left[\sgn(\omega\cdot\epsilon)h\cdot\epsilon\,1_{\{\omega\cdot\epsilon\neq0\}}\right]
			+\E\left[|h\cdot\epsilon|\,1_{\{\omega\cdot\epsilon=0\}}\right] \\
			&=a(\omega)\cdot h
			+\E\left[|h\cdot\epsilon|\,1_{\{\omega\cdot\epsilon=0\}}\right].
		\end{align*}
	\end{proof}

\subsubsection{Proof of Lemma~\ref{lem:equal-blocks}}\label{app:proof-lem-equal-blocks}
\begin{proof}[Proof of Lemma~\ref{lem:equal-blocks}]
		Put $N=2m+2$.  The formulas for $r$ and $s$ give
		\[
		2r^2+2ms^2=1,
		\qquad
		2r=\frac{2-\tau^2}{\sqrt2},
		\]
		so $\omega\in\mathcal C_n(\tau)$.  If $\tau=0$, then
		$\mathcal C_n(0)=\{\omega^*\}$ and the conclusion is immediate.  We may therefore assume that $\tau>0$.  Since $\omega$ is proper, this implies $s>0$.
		
		Define
		\[
		L(h):=\E\left[|h\cdot\epsilon|1_{\{\omega\cdot\epsilon=0\}}\right]
		\]
		and let
		\[
		V:=\left\{h\in\mathbb R^n:
		h_1+h_2=0,\quad \sum_{i=3}^{N}h_i=0\right\}.
		\]
		These are exactly the tangent conditions \eqref{eq:tangent-q} and
		\eqref{eq:tangent-sphere} at $\omega$.
		
		We first show that $L(h)>0$ for every nonzero $h\in V$.  Suppose that
		$h\in V$ and $L(h)=0$.  Then $h\cdot\epsilon=0$ for every sign pattern satisfying
		\[
		\epsilon_1+\epsilon_2=0,
		\qquad
		\sum_{i=3}^{N}\epsilon_i=0,
		\]
		because all these sign patterns belong to the event
		$\{\omega\cdot\epsilon=0\}$.  The signs in the coordinates $i>N$ are arbitrary, so changing one of them gives $h_i=0$ for every $i>N$.  Comparing the two choices $(\epsilon_1,\epsilon_2)=(1,-1)$ and $(-1,1)$ gives $h_1=h_2$.  Similarly, for any $j,k\in\{3,\ldots,N\}$, choose a balanced sign pattern with $\epsilon_j=1$ and $\epsilon_k=-1$ and compare it with the pattern obtained by interchanging these two signs.  This gives $h_j=h_k$.  Thus
		\[
		h_1=h_2,
		\qquad
		h_3=\cdots=h_N,
		\qquad
		h_i=0\quad(i>N).
		\]
		The two defining equations of $V$ now imply $h=0$.  Hence $L$ has trivial kernel on $V$.  Since $V$ is finite-dimensional and $L$ is continuous and positively homogeneous, compactness of the unit sphere in $V$ gives a constant $c>0$ such that
		\begin{equation}\label{eq:equal-block-L-lower}
		L(h)\geq c\|h\|_2
		\qquad\text{for every }h\in V.
		\end{equation}
		
		We now pass from tangent directions to nearby points of
		$\mathcal C_n(\tau)$.  Let $\eta=\omega+v\in\mathcal C_n(\tau)$ be close to $\omega$.  Since the sum of the first two coordinates is fixed on $\mathcal C_n(\tau)$,
		\begin{equation}\label{eq:equal-block-v-q}
		v_1+v_2=0.
		\end{equation}
		Moreover, $\|\omega+v\|_2=\|\omega\|_2=1$ gives
		\[
		2\omega\cdot v+\|v\|_2^2=0.
		\]
		Using \eqref{eq:equal-block-v-q} and the form of $\omega$, we obtain
		\begin{equation}\label{eq:equal-block-v-sum}
		\sum_{i=3}^{N}v_i=-\frac{\|v\|_2^2}{2s}.
		\end{equation}
		Thus a feasible displacement satisfies the tangent equations up to an error of order $\|v\|_2^2$.
		
		To make this precise, let $u\in\mathbb R^n$ have value $1/(2m)$ in the coordinates $3,\ldots,N$ and value $0$ elsewhere, and put
		\[
		\lambda:=\sum_{i=3}^{N}v_i,
		\qquad
		h:=v-\lambda u.
		\]
		Then $h\in V$, and \eqref{eq:equal-block-v-sum} gives
		\begin{equation}\label{eq:equal-block-tangent-error}
		\|v-h\|_2
		=\frac{|\lambda|}{\sqrt{2m}}
		=\frac{\|v\|_2^2}{2s\sqrt{2m}}.
		\end{equation}
		For every $w\in\mathbb R^n$, Cauchy--Schwarz gives
		\[
		L(w)\leq\E|w\cdot\epsilon|\leq\|w\|_2.
		\]
		Therefore, by \eqref{eq:equal-block-L-lower} and
		\eqref{eq:equal-block-tangent-error}, there is a constant $C_1>0$ such that
		\begin{equation}\label{eq:equal-block-L-v}
		L(v)\geq L(h)-L(v-h)\geq c\|h\|_2-\|v-h\|_2
		\geq c\|v\|_2-C_1\|v\|_2^2.
		\end{equation}
		
		Finally, by symmetry and the definition preceding Proposition
		\ref{prop:directional}, there are constants $\alpha,\beta$ such that
		\[
		a(\omega)
		=(\alpha,\alpha,
		\underbrace{\beta,\ldots,\beta}_{2m\text{ coordinates}},
		0,\ldots,0).
		\]
		Hence \eqref{eq:equal-block-v-q} and
		\eqref{eq:equal-block-v-sum} imply
		\begin{equation}\label{eq:equal-block-linear-error}
		a(\omega)\cdot v
		=-\frac{\beta}{2s}\|v\|_2^2.
		\end{equation}
		Because there are only finitely many sign vectors, if $v$ is sufficiently small, then the sign of $(\omega+v)\cdot\epsilon$ agrees with that of $\omega\cdot\epsilon$ whenever $\omega\cdot\epsilon\neq0$.  For such $v$ we have the exact identity
		\[
		K_n(\omega+v)-K_n(\omega)
		=a(\omega)\cdot v+L(v).
		\]
		Combining this identity with \eqref{eq:equal-block-L-v} and
		\eqref{eq:equal-block-linear-error}, we obtain a constant $C_2>0$ such that
		\[
		K_n(\omega+v)-K_n(\omega)
		\geq c\|v\|_2-C_2\|v\|_2^2.
		\]
		The right-hand side is strictly positive for every sufficiently small nonzero $v$.  Therefore $\omega$ is a strict local minimizer of $K_n$ on $\mathcal C_n(\tau)$.
	\end{proof}

\subsubsection{Proof of Lemma~\ref{lem:spike}}\label{app:proof-lem-spike}
\begin{proof}[Proof of Lemma~\ref{lem:spike}]
		The definitions of $a_n(q)$ and $b_n(q)$ imply
		\[
		a_n(q)+b_n(q)=q,
		\qquad
		a_n(q)^2+(n-1)b_n(q)^2=1.
		\]
		Moreover, the condition $q<\sqrt{n/(n-1)}$ gives
		\[
		a_n(q)>(n-1)b_n(q).
		\]
		Hence the first coordinate dominates the sum of the remaining coordinates, and therefore
		\[
		K_n(D_n(q))=a_n(q).
		\]
		
		Now consider a proper vector $\omega\in\mathcal{C}_n(\tau)$ in a sufficiently small neighborhood of $D_n(q)$.  Write $y=\omega_2$.  Since $\omega_1=q-y$ and $\omega_i\leq y$ for $i\geq2$, we have
		\[
		1=\|\omega\|_2^2\leq(q-y)^2+(n-1)y^2.
		\]
		At $y=b_n(q)$ the right-hand side is equal to $1$, and its derivative is negative.  Therefore, in a small neighborhood of $b_n(q)$, the preceding inequality forces
		\[
		y\leq b_n(q).
		\]
		On the other hand, conditioning on $\epsilon_1$ and using Jensen's inequality gives
		\[
		K_n(\omega)\geq\omega_1=q-y.
		\]
		Consequently
		\[
		K_n(\omega)\geq q-y\geq q-b_n(q)=a_n(q)=K_n(D_n(q)).
		\]
		Equality can hold in this neighborhood only when $y=b_n(q)$ and all remaining coordinates also equal $b_n(q)$, namely when $\omega=D_n(q)$.  Thus $D_n(q)$ is a strict local minimizer.
	\end{proof}

\subsubsection{Proof of Proposition~\ref{prop:tau-regularity}}\label{app:proof-prop-tau-regularity}
\begin{proof}[Proof of Proposition~\ref{prop:tau-regularity}]
		Every vector in $\mathbb S^{n-1}_+$ has at least one decreasing
		ordering of its coordinates, so the chambers $\mathfrak C_\pi$ cover
		$\mathbb S^{n-1}_+$.  If all coordinates are distinct, that ordering is
		unique; hence the relative interiors are disjoint.  Permuting coordinates
		gives all $n!$ chambers, and the identity ordering is precisely the proper
		chamber.  Applying coordinatewise sign changes gives the stated
		decomposition of the whole sphere.

		Now suppose that $\omega$ is proper and set
		\[
			q=\omega_1+\omega_2.
		\]
		Since $\omega^*=(1/\sqrt2,1/\sqrt2,0,\ldots,0)$ and both vectors
		have norm one,
		\[
			\tau^2=\|\omega-\omega^*\|_2^2
			=2-\sqrt2\,q,
			\qquad
			q=\frac{2-\tau^2}{\sqrt2}.
		\]
		First suppose that $q\geq1$, or equivalently $0\leq\tau\leq d_0$.
		Because $\sum_{i=3}^n\omega_i^2\geq0$, we have
		$\omega_1^2+\omega_2^2\leq1$.  Therefore
		\[
			(\omega_1-\omega_2)^2
			=2(\omega_1^2+\omega_2^2)-q^2
			\leq2-q^2.
		\]
		Properness gives $\omega_1\geq\omega_2$, and consequently
		\[
			\omega_1
			=\frac{q+(\omega_1-\omega_2)}2
			\leq\frac{q+\sqrt{2-q^2}}2
			=\frac{2-\tau^2+\tau\sqrt{4-\tau^2}}{2\sqrt2}.
		\]
		This bound is attained by the proper two-coordinate unit vector
		\[
			\left(
			\frac{q+\sqrt{2-q^2}}2,
			\frac{q-\sqrt{2-q^2}}2,
			0,\ldots,0
			\right).
		\]

		Next suppose that $0<q\leq1$, equivalently
		$d_0\leq\tau<\sqrt2$.  Since $\omega_2\geq0$,
		\[
			\omega_1\leq\omega_1+\omega_2=q
			=\frac{2-\tau^2}{\sqrt2}.
		\]
		Together with $0\leq\omega_i\leq\omega_1$, this proves the asserted
		$\sigma(\tau)$-regularity on the whole interval.

		It remains to prove that the second branch is dimension-independently
		sharp.  Fix $0<q<1$.  For every sufficiently large integer $m$, put
		\[
			y_m=\frac{q+\sqrt{m-(m-1)q^2}}{m},
			\qquad
			x_m=q-y_m,
		\]
		and define
		\[
			\omega^{(m)}
			=(x_m,\underbrace{y_m,\ldots,y_m}_{m-1\text{ coordinates}})
			\in\mathbb R^m.
		\]
		For sufficiently large $m$ we have $x_m\geq y_m\geq0$.  The defining
		quadratic equation for $y_m$ gives
		\[
			x_m^2+(m-1)y_m^2=1,
			\qquad x_m+y_m=q.
		\]
		Thus $\omega^{(m)}$ is a proper unit vector with
		$\tau(\omega^{(m)})=\tau$.  Moreover,
		\[
			y_m=\frac{\sqrt{1-q^2}}{\sqrt m}+O(m^{-1}),
			\qquad
			x_m=q-\frac{\sqrt{1-q^2}}{\sqrt m}+O(m^{-1}),
		\]
		so $\|\omega^{(m)}\|_\infty=x_m\to q$.  Hence no smaller
		dimension-free bound is possible on the second branch.

		Finally, if $\tau(\omega^{(k)})\to\sqrt2$, then eventually
		$\tau(\omega^{(k)})\geq d_0$, and the second branch gives
		\[
			\|\omega^{(k)}\|_\infty
			\leq\frac{2-\tau(\omega^{(k)})^2}{\sqrt2}
			\longrightarrow0.
		\]
	\end{proof}

\section{Geometric structure behind the phase transitions of $K_n(\omega)$}
\label{app:phase-geometry}

Theorems~\ref{thm:four} and \ref{th:six-local} establish a transition at
$d_{\max}(4)$. We explain why the further even-dimensional endpoints are natural
candidates for analogous transitions. Three features interact: the support restrictions imposed by
the distance constraint, the balanced sign patterns of even uniform
vectors, and the preference for introducing a new pair of coordinates.
The statements below make these features precise. They do not establish
global phase transitions beyond the interval in Theorem~\ref{th:six-local}.

\subsection{Distance slices and support thresholds}
For $m\geq4$, write
\[
 u_m=(\underbrace{m^{-1/2},\ldots,m^{-1/2}}_{m},0,\ldots,0),
 \qquad d_m=d_{\max}(m)=\sqrt{2-2\sqrt{2/m}}.
\]
The vector $u_m$ may be embedded in any ambient dimension $n\geq m$.
For a proper unit vector, the distance constraint is equivalent to fixing
\begin{equation}\label{eq:geometry-q-slice}
 q=\omega_1+\omega_2=\frac{2-\tau^2}{\sqrt2}.
\end{equation}
Thus increasing $\tau$ moves the hyperplane $\omega_1+\omega_2=q$
through the proper part of the unit sphere, with $q$ decreasing.

\begin{geometryproposition}[Support thresholds]\label{prop:geometry-support}
Let $\omega$ be a proper unit vector with at most $m$ nonzero coordinates,
where $4\leq m\leq n$. Then $q\geq2/\sqrt m$ and $\tau(\omega)\leq d_m$.
If $m>4$, equality holds if and only if $\omega=u_m$.
For $m=4$, equality also occurs at $(1,0,\ldots,0)$.
Consequently, $\tau>d_m$ forces more than $m$ nonzero coordinates.
\end{geometryproposition}
\begin{proof}
Put $y=\omega_2\in[0,q/2]$ and use properness to obtain
\[
 1\leq(q-y)^2+(m-1)y^2.
\]
The right side is strictly convex in $y$, so its maximum on this interval
is $\max\{q^2,mq^2/4\}=mq^2/4$. This proves the claimed bound on $q$,
and \eqref{eq:geometry-q-slice} gives the bound on $\tau$.
When $m>4$ and $q=2/\sqrt m$, only $y=q/2$ can attain the maximum.
Equality in the norm estimate then forces all $m$ coordinates to be
$1/\sqrt m$. At $m=4$, the other endpoint $y=0$ gives the dictator;
these are the only equality cases.
\end{proof}

The support thresholds are therefore genuine changes in the geometry of
the feasible set. Nevertheless, this observation alone does not explain
why even dimensions should be distinguished: it also applies to odd $m$.
The missing ingredient is the parity of the zero sums.

\subsection{Balanced sign patterns distinguish even dimensions}
Let $S_m=\epsilon_1+\cdots+\epsilon_m$. For even $m$, put
\begin{equation}\label{eq:geometry-balance-mass}
 p_m=\mathbb P(S_m=0)=2^{-m}\binom{m}{m/2},
 \qquad \E|S_m|=mp_m,\qquad K_n(u_m)=\sqrt m\,p_m.
\end{equation}
In coefficient space, $K_n$ is linear within every region cut out by the
hyperplanes $\omega\cdot\epsilon=0$. The configuration $u_m$ lies on
many such hyperplanes when $m$ is even, whereas it lies on none when
$m$ is odd. The following local distinction is a consequence.

\begin{geometryproposition}[Parity and local extremality]
\label{prop:geometry-parity}
Fix $n\geq m\geq4$. If $m$ is even, $u_m$ is a strict local minimizer of
$K_n$ on the unit sphere, and hence on $\mathcal C_n(d_m)$.
If $m$ is odd, $u_m$ is a strict local maximizer on the unit sphere.
Moreover, when $m\geq5$ is odd and $n\geq m+1$, it is not a local
minimizer on $\mathcal C_n(d_m)$.
\end{geometryproposition}
\begin{proof}
Use $\sgn(0)=0$ in this proof, and let
$a=\E[\epsilon\,\sgn(u_m\cdot\epsilon)]$.
Symmetry gives $a=K_n(u_m)u_m$.
Since there are finitely many sign patterns, all nonzero sums keep their
sign in a sufficiently small neighborhood of $u_m$. For such $v$,
\begin{equation}\label{eq:geometry-local-expansion}
 K_n(u_m+v)-K_n(u_m)
 =K_n(u_m)\langle u_m,v\rangle+L_m(v),\qquad
 L_m(v)=\E\bigl[|v\cdot\epsilon|1_{\{S_m=0\}}\bigr].
\end{equation}
For unit $u_m+v$, $\langle u_m,v\rangle=-\|v\|_2^2/2$.
If $m$ is odd, $L_m=0$, so the difference is strictly negative for
nonzero $v$.

Suppose $m$ is even. On the tangent space $u_m^\perp$, the function
$L_m$ vanishes only at zero. Indeed, $L_m(h)=0$ implies
$h\cdot\epsilon=0$ for every balanced head pattern and every tail pattern.
Changing a single tail sign gives $h_i=0$ for $i>m$. Swapping opposite
signs in a balanced head pattern gives $h_i=h_j$ for $i,j\leq m$.
Tangency then forces these common entries to be zero.
By continuity and compactness, $L_m(h)\geq c_{m,n}\|h\|_2$ on
$u_m^\perp$, for some $c_{m,n}>0$.
Since $L_m(v)=L_m(v-\langle u_m,v\rangle u_m)$, the positive term in
\eqref{eq:geometry-local-expansion} is of order $\|v\|_2$, while the
negative term is of order $\|v\|_2^2$. This proves strict local minimality.

To check that odd $m$ also allows a decreasing path on the fixed-distance
slice, put $r=1/\sqrt m$ and, for $0<t\leq r/\sqrt2$, define
\[
 w(t)=(\underbrace{r,\ldots,r}_{m-1},\sqrt{r^2-t^2},t,0,\ldots,0).
\]
These vectors are proper, unit, have $\omega_1+\omega_2=2r$, and tend
to $u_m$. Put $p=2^{-(m-1)}\binom{m-1}{(m-1)/2}$.
Every nonzero sum of the first $m-1$ coordinates has magnitude at least
$2r$, whereas the last two coordinates contribute at most
$\sqrt{r^2-t^2}+t\leq\sqrt2r<2r$.
Conditioning on the even head therefore gives
\[
 K_n(w(t))=p\bigl((m-1)r+\sqrt{r^2-t^2}\bigr)
 <pmr=K_n(u_m).
\]
Thus $u_m$ is not locally minimizing on this slice.
\end{proof}

This parity distinction gives a geometric reason to expect low-dimensional
minimizing branches to pass through even uniform configurations.
Odd support thresholds do not have the same local-minimum mechanism when
additional coordinates are available.

\subsection{Pair insertion and the square-root scale}
The balanced sign patterns also explain why two new coordinates are
favored when one leaves an even uniform configuration.

\begin{geometryproposition}[The exact profile with an equal head]
\label{prop:geometry-pair-insertion}
Let $m\geq4$ be even, $n\geq m+2$, and $d_m\leq\tau\leq d_{m+2}$.
Put
\[
 q=\frac{2-\tau^2}{\sqrt2},\qquad r=\frac q2,\qquad
 t=\sqrt{1-mr^2},\qquad s=\frac{t}{\sqrt2}.
\]
Among the vectors in $\mathcal C_n(\tau)$ whose first $m$ coordinates
are equal, the unique minimizer is
\[
 v_m(\tau)=(\underbrace{r,\ldots,r}_{m},s,s,0,\ldots,0),
\]
and its value is
\begin{equation}\label{eq:geometry-pair-value}
 K_n(v_m(\tau))=p_m\left(\frac{mq}{2}
                 +\frac1{\sqrt2}\sqrt{1-\frac{mq^2}{4}}\right).
\end{equation}
In particular, as $\varepsilon\downarrow0$,
\begin{equation}\label{eq:geometry-pair-cusp}
 K_n(v_m(d_m+\varepsilon))
 =K_n(u_m)+\beta_m\sqrt\varepsilon+O(\varepsilon),\qquad
 \beta_m=p_m\left(\frac m2\right)^{1/4}\sqrt{d_m}>0.
\end{equation}
\end{geometryproposition}
\begin{proof}
Every vector in the stated class is $(r,\ldots,r,y)$ with
$\|y\|_2=t$. Let $Y$ be its tail Rademacher sum. Conditional on a
nonzero value of $S_m$, Jensen's inequality gives
$\E_Y|rS_m+Y|\geq r|S_m|$; on $S_m=0$, the contribution is $K(y)$.
Using \eqref{eq:geometry-balance-mass} and the sharp Khintchine inequality,
\[
 K_n(r,\ldots,r,y)\geq r\E|S_m|+p_mK(y)
 \geq p_m\left(mr+\frac t{\sqrt2}\right).
\]
The interval for $\tau$ is exactly
$2/\sqrt{m+2}\leq q\leq2/\sqrt m$, which ensures $s\leq r$.
For the tail $(s,s)$, its maximum magnitude is $2s\leq2r$;
hence it cannot change the sign of a nonzero head sum. All inequalities
are then equalities, giving \eqref{eq:geometry-pair-value}.
For $t>0$, equality in the tail Khintchine inequality forces a pair of
equal nonzero tail coordinates. Properness fixes their positions.
At $t=0$, the only vector in the class is $u_m$.

At $\tau=d_m$, $q=2/\sqrt m$ and $t=0$. Differentiation gives
\[
 t^2=1-\frac{mq^2}{4}
 =\sqrt{2m}\,d_m\varepsilon+O(\varepsilon^2),\qquad
 \frac{mq}{2}=\sqrt m+O(\varepsilon).
\]
Substitution into \eqref{eq:geometry-pair-value} proves
\eqref{eq:geometry-pair-cusp}.
\end{proof}

The square-root scale thus has a direct source. The distance constraint
creates tail \emph{variance} proportional to $\tau-d_m$, whereas the
balanced head patterns contribute the tail's expected absolute value,
which is proportional to its standard deviation. Among all tails, the
sharp Khintchine inequality favors a balanced pair. A single new
coordinate has $K(y)=\|y\|_2$, while a balanced pair has
$K(y)=\|y\|_2/\sqrt2$.

For $m=4$, formula~\eqref{eq:geometry-pair-value} is precisely
\[
 \frac{3q}{4}+\frac{3}{8\sqrt2}\sqrt{1-q^2},
\]
and $\beta_4=\frac38(4-2\sqrt2)^{1/4}$, agreeing with
Corollary~\ref{cor:cusp}. Theorem~\ref{th:six-local} proves that this
restricted minimum is also the global minimum on its explicit initial
interval. For higher even $m$, Proposition~\ref{prop:geometry-pair-insertion}
provides an exact restricted minimum and therefore an \emph{upper bound}
on the unrestricted minimum $M_n(\tau)$.

\subsection{What the proposed higher transitions still require}
The first even uniform configurations give
\[
\begin{array}{c|c|c}
 m&d_m&K_n(u_m)\\ \hline
 4&\sqrt{2-\sqrt2}&3/4\\
 6&\sqrt{2-2/\sqrt3}&15/(8\sqrt6)\\
 8&1&35/(32\sqrt2)
\end{array}
\]
At each $d_m$ with $m$ even, the previous support size reaches its
geometric limit, the uniform point is locally minimizing, and an
optimal equal-head continuation adds a new pair with square-root cost.
Together these facts motivate the sequence of candidate transitions
$d_{\max}(4),d_{\max}(6),d_{\max}(8),\ldots$.

To promote this picture to a theorem about the global profile, one must
also prove that $u_m$ is a global minimizer at $d_m$ and control all
competing configurations for nearby distances, including those with
unequal heads. Local minimality and support restrictions do not supply
these comparisons. Nor do the arguments above exclude additional branch
changes inside the intervals $(d_m,d_{m+2})$: Conjecture~\ref{conj:six}
already predicts an internal change of the six-coordinate family.
The present paper establishes the global transition at $d_4$; the
corresponding global statements at $d_6,d_8,\ldots$ remain conjectural.

\section{Additional properties}\label{app:additional}
\subsection{Local minimizers of $K_n(\omega)$ on $\mathcal{C}_n(\tau)$}
	
	We begin by recalling a structural lemma for $K_n$ on the proper chamber.
	
	\begin{lemma}[{\cite[Lemma 8]{kalarickal2025convexity}}]\label{lem:convex}
		Let
		\[
		T_n=\{\omega\in\mathbb{R}^n:\ \omega_1\geq\omega_2\geq\cdots\geq\omega_n\geq0\}.
		\]
		For every $n\geq1$, there exists a finite set $A_n\subseteq T_n$ such that, for all $\omega\in T_n$,
		\[
		K_n(\omega)=\max_{a\in A_n}a\cdot\omega.
		\]
	\end{lemma}
	
	It is also mentioned in \cite{kalarickal2025convexity} that, if
	\[
	a_j(\omega)=\E\left[\epsilon_j\sgn\left(\sum_{i=1}^n\omega_i\epsilon_i\right)\right],
	\qquad 1\leq j\leq n,
	\]
	then $a(\omega)\in T_n$ and $K_n(\omega)=a(\omega)\cdot\omega$.  Here we use the convention $\sgn(0)=0$.
	
	For a nonzero direction $h\in\mathbb{R}^n$, the one-sided directional derivative is defined as
	\[
	DK_n(\omega;h) := \lim_{\delta\to0^+}
	\frac{K_n(\omega+\delta h)-K_n(\omega)}{\delta}. 
	\]
	
	\begin{proposition}\label{prop:directional}
		For every direction $h\in\mathbb{R}^n$,
		\[
		DK_n(\omega;h)
		=a(\omega)\cdot h+
		\E\left[|h\cdot\epsilon|\,1_{\{\omega\cdot\epsilon=0\}}\right],
		\]
		where $\epsilon=(\epsilon_1,\ldots,\epsilon_n)$.
	\end{proposition}
	
	The proof is given in Appendix~\ref{app:proof-prop-directional}.

	On $\mathcal{C}_n(\tau)$ the quantity $\omega_1+\omega_2$ is fixed, because
	\[
	\omega_1+\omega_2=\frac{2-\tau^2}{\sqrt2}.
	\]
	Therefore, when studying minimizers on $\mathcal{C}_n(\tau)$, we only need to consider first-order directions $h$ with
	\begin{equation}\label{eq:tangent-q}
		h_1+h_2=0.
	\end{equation}
	The normalization $\|\omega\|_2=1$ further gives
	\begin{equation}\label{eq:tangent-sphere}
		\omega\cdot h=\sum_{i=1}^n\omega_i h_i=0.
	\end{equation}
	
	Equations \eqref{eq:tangent-q} and \eqref{eq:tangent-sphere} define the tangent space to the equality constraints. A feasible first-order direction must also respect the tangent cone of the proper chamber. At a local minimizer, $DK_n(\omega;h)\geq0$ for every direction realized by a differentiable feasible curve. In Proposition \ref{prop:directional}, the second term is always nonnegative.  The first term, however, changes sign when $h$ is replaced by $-h$.  A useful way to construct local minimizers is therefore to arrange that $a(\omega)\cdot h=0$ for all tangent directions and then prove strict positivity of the second term.
	
	\begin{lemma}\label{lem:equal-blocks}
		Let $m\in\mathbb{N}_+$ and assume $n\geq 2m+2$.  Let
		\[
		\omega=(r,r,\underbrace{s,\ldots,s}_{2m\text{ coordinates}},0,\ldots,0),
		\]
		where
		\[
		r=\frac{2-\tau^2}{2\sqrt2},
		\qquad
		s=\sqrt{\frac{4\tau^2-\tau^4}{8m}}.
		\]
		Whenever this vector is proper, it is a strict local minimizer of $K_n$ on $\mathcal{C}_n(\tau)$.
	\end{lemma}
	
	The proof is given in Appendix~\ref{app:proof-lem-equal-blocks}.

	This lemma explains why the second family in Conjecture \ref{conj:six} is a natural candidate: it is one of the equal-block local minimizers produced by the directional-derivative criterion.
	There is another type of local minimizer near $\tau=d_{\max}(4)$.
	
	\begin{lemma}\label{lem:spike}
		Let
		\[
		q=\frac{2-\tau^2}{\sqrt2}.
		\]
		For $1\leq q<\sqrt{n/(n-1)}$, define
		\[
		b_n(q)=\frac{q-\sqrt{n-(n-1)q^2}}{n},
		\qquad
		a_n(q)=q-b_n(q).
		\]
		Then
		\[
		D_n(q)=(a_n(q),b_n(q),\ldots,b_n(q))
		\]
		is a strict local minimizer of $K_n$ on $\mathcal{C}_n(\tau)$.
	\end{lemma}
	
	The proof is given in Appendix~\ref{app:proof-lem-spike}.

\subsection{Order chambers and regularity}

	\begin{definition}\label{def:regular}
		For $\sigma>0$, a proper unit vector $\omega\in\mathbb R^n$ is called
		\emph{$\sigma$-regular} if $|\omega_i|\leq\sigma$ for every
		$1\leq i\leq n$.
	\end{definition}

	\begin{proposition}[Order chambers and regularity]\label{prop:tau-regularity}
		Let
		\[
			\mathbb S^{n-1}_+
			=\mathbb S^{n-1}\cap[0,\infty)^n.
		\]
		The coordinate orderings divide $\mathbb S^{n-1}_+$ into the $n!$
		closed chambers
		\[
			\mathfrak C_\pi
			=\bigl\{\omega\in\mathbb S^{n-1}_+:
			\omega_{\pi(1)}\geq\omega_{\pi(2)}\geq\cdots
			\geq\omega_{\pi(n)}\bigr\},
			\qquad \pi\in S_n,
		\]
		whose relative interiors are pairwise disjoint.  The proper chamber is
		$\mathfrak C_{\mathrm{id}}$.  Equivalently, modulo independent sign
		changes, these are the $n!$ coordinate-order chambers on
		$\mathbb S^{n-1}$; if signs are retained, the full signed decomposition
		consists of $2^n n!$ chambers, up to their common boundaries.

		Moreover, if $n\geq2$, $\omega$ is a proper unit vector, and
		$\tau=\tau(\omega)$, then $\omega$ is $\sigma(\tau)$-regular, where
		$d_0=\sqrt{2-\sqrt2}$ and
		\begin{equation}\label{eq:sigma-tau}
			\sigma(\tau)=
			\begin{cases}
			\displaystyle
			\frac{2-\tau^2+\tau\sqrt{4-\tau^2}}{2\sqrt2},
			&0\leq\tau\leq d_0,\\[3mm]
			\displaystyle
			\frac{2-\tau^2}{\sqrt2},
			&d_0\leq\tau<\sqrt2.
			\end{cases}
		\end{equation}
		This is the smallest possible bound that is independent of the dimension.
		In particular, if $\omega^{(k)}$ is any sequence of proper unit vectors,
		possibly of varying dimensions, such that
		$\tau(\omega^{(k)})\to\sqrt2$, then
		\[
			\|\omega^{(k)}\|_\infty\longrightarrow0.
		\]
	\end{proposition}

	The proof is given in Appendix~\ref{app:proof-prop-tau-regularity}.

\nocite{tomaszewski1987simple,figiel1997extremal,boppana2017tomaszewski,chasapis2021ball,havrilla2021sharp,havrilla2023khinchin,nayar2023extremal,melbourne2022quantitative}
\begingroup
\small
\bibliographystyle{plain}
\bibliography{bmyref}
\endgroup
\end{document}